\documentclass[11pt]{article}
\usepackage{amsmath,mathrsfs,amssymb,amsthm,amscd}
\usepackage[]{fontenc}
\usepackage[all]{xy}
\usepackage{hyperref}
\usepackage{graphics}
\usepackage{a4wide}
\usepackage[dvips]{graphicx,color}
\usepackage{setspace}

\newcounter{EQNR}
\renewcommand{\contentsname}{Contents\\
{\footnotesize\normalfont(A table of contents should normally not
be included)}}

\newtheorem{thm}{\bf Theorem}[section]
\newtheorem{Def}[thm]{\bf Definition}
\newtheorem{Lem}[thm]{\bf Lemma}
\newtheorem{Pro}[thm]{\bf Proposition}

\newcommand{\Del}{{\Delta}}
\newcommand{\R}{\mathbb{R}}
\newcommand{\N}{\mathbb{N}}
\newcommand{\s}{\mathbb{S}}
\newcommand{\D}{\mathbb{D}}
\newcommand{\tphi}{\tilde{\phi}}
\newcommand{\del}{{\partial}}
\newcommand{\Vol}{{\text{Vol}}}

\newcommand{\Ric}{{\text{Ric}}}
\newcommand{\diam}{{\text {diam}}}
\newcommand{\skel}{{\text {skel}}}
\newcommand{\cut}{{\text {cut}}}
\newcommand{\Ann}{{\text {Ann}}}
\newcommand{\im}{{\text {Im}}}

\newcommand{\clos}{{\text {clos}}}

\begin{document}

\title{Volume growth and the topology of manifolds with nonnegative Ricci
curvature}
\author{Michael Munn}
\date{}

\maketitle

\begin{abstract}\noindent
Let $M^n$ be a complete, open Riemannian manifold with $\Ric \geq 0$. In
1994, Grigori Perelman showed that there exists a constant
$\delta_{n}>0$, depending only on the dimension of the manifold,
such that if the volume growth satisfies $\alpha_M := \lim_{r \rightarrow \infty}
\frac{\Vol(B_p(r))}{\omega_n r^n} \geq 1-\delta_{n}$, then $M^n$ is
contractible. Here we employ the techniques of Perelman to find
specific lower bounds for the volume growth, $\alpha(k,n)$, depending only on $k$ and $n$, which guarantee the individual $k$-homotopy group of $M^n$ is trivial.
\end{abstract}

\section{Introduction}\label{Chapter-Intro}

\doublespacing 

Let $M^n$ be an $n$-dimensional complete Riemannian manifold with
nonnegative Ricci curvature. For a base point $p \in M^n$, denote by
$B_p(r)$ the open geodesic ball in $M^n$ centered at $p$ and with
radius $r$. Let $\Vol(B_p(r))$ denote the volume of $B_p(r)$ and
denote by $\omega_n$ the volume of the unit ball in Euclidean space.
By the Bishop-Gromov Relative Volume Comparison Theorem, \cite{BC,
GLP}, the function $r \hspace{-.1cm} \rightarrow \Vol(B_p(r))/\omega_n r^n$ is
non-increasing and bounded above by 1.

\begin{Def}
Define $\alpha_M$, the volume growth of $M^n$, as
$$ \alpha_M := \lim_{r \rightarrow \infty} \frac{\Vol(B_p(r))}{\omega_n r^n}.$$

The manifold $M^n$ is said to have Euclidean (or large) volume growth when $\alpha_M > 0$.
\end{Def}

The constant $\alpha_M$ is a global geometric invariant of $M^n$,
i.e. it is independent of base point. Also, when $\alpha_M  > 0$,
$$
\Vol(B_p(r) \geq \alpha_M \omega_n r^n, \qquad \text{for all } p
\in M ~\text{and for all } r>0.
$$
It follows from the Bishop-Gromov Volume Comparison
Theorem \cite{BC, GLP} that $\alpha_M = 1$ implies $M^n$ is isometric to
$\R^n$. 

In this paper, we study complete manifolds with $\Ric_M \geq 0$ and
$\alpha_M >0$. Anderson ~\cite{A} and Li ~\cite{Li} have
independently shown that the order of $\pi_1(M^n)$ is bounded from
above by $\frac{1}{\alpha_M}$. In particular, if $\alpha_M >
\frac{1}{2}$, then $\pi_1(M^n) =0$. Furthermore, Zhu \cite{Z1} has
shown that when $n=3$, if $\alpha_M >0$, then $M^3$ is contractible.
It is interesting to note that this is not the case when $n=4$ as
Menguy \cite{M} has constructed examples of 4-manifolds with large
volume growth and infinite topological type based on an example by
Perelman \cite{P1}. In 1994, Perelman ~\cite{P2} proved that there
exists a small constant $\delta_n>0$ which depends only on the
dimension $n \geq 2$ of the manifold, such that if $\alpha_M \geq 1-
\delta_n$, then $M^n$ is contractible. It was later shown by Cheeger
and Colding \cite{ChCoI} that the conditions in Perelman's theorem
are enough to show that $M^n$ is $C^{1, \alpha}$ diffeomorphic to
$\R^n$. In this paper, we follow the method of proof in Perelman's
theorem. Employing this method, we determine specific bounds on
$\alpha_M$ which imply the individual $k$-th homotopy groups of the
manifold are trivial.

We prove

\begin{thm}\label{Theorem-Main}\textit{
Let $M^n$ be a complete Riemannian manifold with $\Ric \geq 0$. If
$$\alpha_M > \alpha(k,n),$$
where $\alpha(k,n)$ are the constants given in Table \ref{Table of alphas},
then $\pi_k(M^n)=0$.}
\end{thm}

{\bf Remark.}
Table \ref{Table of alphas} contains values of $\alpha(k,n)$ for $1 \leq k \leq 3$ and $1 \leq n \leq 10$. In general, the value of $\alpha(k,n)$ is determined by Equation (\ref{alphas}), where the function $h_{k,n}(x)$ and the values of $\delta_{k,n}$ are defined in Definition \ref{Def-six* function}.

In section \ref{Section-Background}, we state general
results from Riemannian geometry that will be required for the
proof. The key ingredients are the excess estimate of
Abresch-Gromoll, the Bishop-Gromov Volume Comparison Theorem, and a
Maximal Volume Lemma of Perelman [Lemma \ref{Lemma-Perelman's
Maximal Volume}]. 

In section \ref{Chapter-almost equicontinuity}, we
apply the theory of almost equicontinuity from \cite{So} to prove a
general Homotopy Construction Theorem [Theorem
\ref{Theorem-HomotopyConstruction}] that will be needed when
constructing the homotopies for Theorem \ref{Theorem-Main}.

In section \ref{Section-General Argument}, we prove Theorem
\ref{Theorem-Main} using a double induction argument for the general
case. This argument follows Perelman's except that we carefully
determine the necessary constants to build each step. Perelman's double induction argument is built from two lemmas each
of which depends on a parameter $k\in \N$.  The Main Lemma($k$) [Lemma
\ref{Lemma-Main}] says that given a constant $c>1$ and an appropriate estimate on volume
growth, any given continuous function $f: \s^k \to B_p(R)$ can be
extended to a continuous function $g: \D^{k+1} \to B_p(cR)$ .  This lemma is proven by defining intermediate
functions $g_j$ on finer and finer nets in $\D^{k+1}$. To define $g_j$ on these
nets one uses the Moving In Lemma, described below. To prove the limit $g(x) =
\lim_{j \to \infty} g_j(x)$ exists and is continuous, we apply results from section
\ref{Chapter-almost equicontinuity}.

The Moving In Lemma($k$) [Lemma \ref{Lemma-Moving In}] states that given a constant $d_0 >0$ and a
map $\phi: \s^k \to B_q(\rho)$ then with an appropriate bound on
volume growth one can move $\phi$ inward obtaining a new map $\tphi:
\s^k \to B_q((1-d_0)\rho)$. The new map $\tphi$ is uniformly close
to the map $\phi$ with respect to the radius $\rho$. The maps $\phi$
and $\tphi$ are not necessarily homotopic; however, a homotopy is
constructed by controlling precisely the uniform closeness of these
maps on smaller and smaller scales. The Moving In Lemma($k$) and
Main Lemma($i$), for $i=0,..,k-1$, are used to produce finer and
finer nets that then converge on the homotopy required for Main
Lemma($k$). Moving In Lemma($k$) is proven by constructing the map
$\tphi$ inductively on successive $i$-skeleta of a triangulation of
$\s^k$. The conclusion of Main Lemma($i$), for $i=0,..,k-1$, is
needed in the induction step of the proof of Moving In Lemma($k$).

The key place in the argument where the volume growth bound is
introduced occurs in the proof of the Moving In Lemma; specifically,
in producing a small, thin triangle in an advantageous location.
However, due to the double inductive argument, and the fact that
lower dimensional lemmas are applied on a variety of scales where
the choice of $c$ and $d_0$ depend on $n$ and $k$, the actual
estimate on the volume is produced using inductively defined
functions $\beta(k, c, n)$ [Definition \ref{Def-beta volume growth}]
and constants $C_{k,n}$ [Definition \ref{Def-C(k,n) constants}].

In the appendix, we complete our analysis of $\beta(k, c, n)$ to find the optimal bounds, $\alpha(k,n)$, over all constants $c>1$. Through this analysis we are able to construct a table of values containing the optimal lower bounds for the volume growth, as stated in Theorem \ref{Theorem-Main}, which guarantee the $k$-th homotopy group is trivial. The bounds that we obtain are the best that can be achieved via Perelman's method. A portion of this analysis was done using Mathematica 6. The code for these commands is available in \cite{Mu}.

\vskip .10in \textbf{Acknowledgements.} I would like to thank Professors Christina Sormani and Isaac Chavel for their advice, encouragement, and insight while working on this project. 

\subsection{Background}\label{Section-Background}
Here we review two facts from the Riemannian geometry of manifolds
with non-negative Ricci curvature. Let $M^n$ be a complete
Riemannian manifold with $\Ric \geq 0$.

\begin{thm}\label{Theorem-Abresch-Gromoll Excess}
{\em [Abresch-Gromoll Excess Theorem]}. \textit{ Let $p, q \in M^n$ and
let $\overline{pq}$ be a minimal geodesic connecting $p$ and $q$. For any
$x \in M^n$, we define the excess function with respect to $p$ and
$q$ as
$$e_{p,q}(x) = d(p,x) + d(q, x) - d(p,q).$$
Define $h(x) = d(x,\overline{pq})$ and set $s(x) = \min \left\{d(p,x), d(q,x) \right\}$.
If $h(x) \leq s(x) / 2$, then
$$e_{p,q}(x)\leq 8 \left( \frac{h(x)^n}{s(x)} \right)^{1/{n-1}} = 8
\left( \frac{h(x)}{s(x)} \right)^{1/{n-1}} h(x).$$}
\end{thm}

This excess estimate is due to Abresch-Gromoll \cite{AbGr} (c.f.
\cite{Ch}).

\begin{Def}\label{Def-gamma}
For constants $c>1$, $\epsilon > 0$ and $n \in \N$, define
$$\gamma(c,\epsilon,n) = \left[1 + \left(\frac{c}{\epsilon}\right)^n \right]^{-1}.$$
\end{Def}

\begin{Lem}\label{Lemma-Perelman's Maximal Volume}
{\em [Perelman's Maximal Volume Lemma]}. Let $p \in M^n$, $R>0$,
for any constants $c_1>1$ and $\epsilon>0$, if $\alpha_M > 1 -
\gamma(c_1, \epsilon, n)$, then for every $a \in B_p(R)$, there
exist $q \in M^n \setminus B_p(c_1 R)$ such that $d(a,\overline{pq})
\leq \epsilon R$, where $\overline{pq}$ denotes a minimal geodesic
connecting $p$ and $q$.
\end{Lem}

This fact was observed without proof by Perelman in \cite{P2}. Our statement and proof differ in that we utilize the global volume growth control on $\alpha_M$ rather than only a local volume bound in a neighborhood of $B_p(c_1 R)$. This global bound allows us to determine an expression for $\gamma$ not given in \cite{P2}. The proof of Perelman's original statement follows from the proof of the Bishop-Gromov Volume Comparison Theorem and can be found in \cite{Z2}.

\begin{proof}
Let $c_2 > c_1 >1$ be finite constants. Define $\Gamma \equiv
\{\dot{\sigma} |~ d(a, \sigma) \leq \epsilon R\} \subset
\s_p^{n-1}(M^n) \subset T_pM^n$, where $\sigma$ denotes a minimal geodesic in $M^n$ and $\sigma$ its velocity vector. Suppose that for all $v \in \Gamma$, we have
$\cut(v) < c_1 R$. In what follows, we determine an upper bound on the volume growth, $\alpha_M$, which would allow such a contradiction to occur. In turn, by requesting the volume growth be bounded below by this upper bound, the lemma will follow.

By definition, we have
\begin{eqnarray*}
\Vol(B_p(c_2R)) &=& \int_{\Gamma} \int_0^{\min\{\cut(v), c_2R \}}
A_{M^n}(t, v) dt dv \\
&& \hspace{3.8cm} + \int_{\s^{n-1}\setminus \Gamma}
\int_0^{\min\{\cut(v),c_2R\}}A_{M^n}(t,v) dtdv\\
&\leq& \Vol(\Gamma) \int_{0}^{c_1R} A^0(t)dt +
\Vol(\s^{n-1}\setminus
\Gamma) \int_{0}^{c_2R}A^0(t)dt\\
&=& \Vol(\s^{n-1})\int_{0}^{c_2R}A^0(t)dt - \Vol(\Gamma) \left( \left(\int_{0}^{c_2R} - \int_{0}^{c_1R} \right) A^0(t)dt \right)\\
&=& -\Vol(\Gamma) \int_{c_1R}^{c_2R}A^0(t)dt + \Vol(\s^{n-1}) \int_{0}^{c_2R}A^0(t)dt\\
&=& -\Vol(\Gamma) \int_{c_1R}^{c_2R}A^0(t)dt + \Vol(B^0(c_2R)).
\end{eqnarray*}
Here $A_{M^n}(t,v)$ denotes the volume element on $M^n$ and $A^0(t)$
denotes the volume element on $\R^n$; that is, $A^0(t) =t^{n-1}$. From
the assumption on the volume growth, we have that $\Vol(B_p(c_2R)) \geq (1-\gamma)
\Vol(B^0(c_2R))$ and therefore
\begin{eqnarray}
(1-\gamma)\Vol(B^0(c_2R)) &\leq& -\Vol(\Gamma) \int_{c_1R}^{c_2R}
A^0(t)
+ \Vol(B^0(c_2R))\\
\gamma \Vol(B^0(c_2R)) &\geq& \Vol(\Gamma) \int_{c_1R}^{c_2R}
A^0(t)dt\\
\Vol(\Gamma) &\leq& \gamma
\frac{\Vol(B^0(c_2R))}{\int_{c_1R}^{c_2R}A^0(t)dt}.
\end{eqnarray}
On the other hand, since $B_a(\epsilon R) \subset \Ann_{\Gamma}(p;
0, c_1R)$, it follows that $\Vol(B_a(\epsilon R)) \leq \Vol(\Gamma)
\int_{0}^{c_1R} A^0(t)dt$. Hence
\begin{equation}\label{upper bound for B_a(epsilonR)}
\Vol(B_a(\epsilon R)) \leq \gamma \Vol(B^0(c_2R))
\frac{\int_{0}^{c_1R} A^0(t)dt}{\int_{c_1R}^{c_2R} A^0(t)dt}.
\end{equation}

Furthermore, since $B_p(c_2R) \subset B_a(R + c_2R)$, we know that
$$\frac{\Vol(B_p(c_2R))}{\Vol(B_a(\epsilon R)} \leq \frac{\Vol(B_a(R+c_2R))}{\Vol(B_a(\epsilon R))} \leq \frac{(R+c_2R)^n}{(\epsilon R)^n};$$
and therefore,
\begin{eqnarray}
\Vol(B_a(\epsilon R)) &\geq& \Vol(B_a(R+c_2R))
\frac{(\epsilon R)^n}{(R+c_2R)^n}\\
&\geq& \Vol(B_p(c_2R))
\frac{\epsilon^n}{(1+c_2)^n}\\
&\geq& (1-\gamma) \Vol(B^0(c_2R))
\frac{\epsilon^n}{(1+c_2)^n}.\label{lower bound for B_a(epsilonR)}
\end{eqnarray}
Combining (\ref{upper bound for B_a(epsilonR)}) and (\ref{lower bound
for B_a(epsilonR)}), we get
\begin{eqnarray}
(1-\gamma) \Vol(B^0(c_2R)) \frac{\epsilon^n}{(1+c_2)^n} &\leq&
\gamma \Vol(B^0(c_2R)) \frac{\int_{0}^{c_1R}
A^0(t)dt}{\int_{c_1R}^{c_2R} A^0(t)dt}\\
\left(\frac{\epsilon}{1+c_2} \right)^n - \gamma
\left(\frac{\epsilon}{1+c_2} \right)^n &\leq& \gamma
\frac{c_1^n}{c_2^n - c_1^n}\\
\left(\frac{\epsilon}{1+c_2} \right)^n &\leq& \gamma
\left[\frac{c_1^n}{c_2^n - c_1^n} +
\left(\frac{\epsilon}{1+c_2}\right)^n \right]. \label{gamma bound}
\end{eqnarray}

By solving (\ref{gamma bound}) for $\gamma$, we can deduce a lower bound for $\gamma$ dependent only on the constants $c_2, c_1, \epsilon$ and $n$. That is,
\begin{eqnarray}
\gamma &\geq&
\left(\frac{\epsilon}{1+c_2}\right)^n\left[\frac{c_1^n}{c_2^n -
c_1^n} + \left(\frac{\epsilon}{1+c_2}\right)^n \right]^{-1}\\
&=& \left[1+ \frac{c_1^n}{c_2^n - c_1^n}
\left(\frac{1+c_2}{\epsilon}\right)^n \right]^{-1}.
\end{eqnarray}

Note that, throughout the proof we required a restriction on the volume growth
only within the larger ball $B_p(c_2R)$. Since $\alpha_M$ is a
global restriction on volume growth, it is possible to take $c_2
\rightarrow \infty$ and thus refine the lower bound on $\gamma$ determined above. Since
$$\lim_{c_2 \rightarrow \infty} \left[1+ \frac{c_1^n}{c_2^n - c_1^n}
\left(\frac{1+c_2}{\epsilon}\right)^n \right]^{-1} = \left[1+
\frac{c_1^n}{\epsilon^n}\right]^{-1},$$ the above lower bound on $\gamma$ can be expressed more simply as
$$
\gamma \geq \left[1+
\frac{c_1^n}{\epsilon^n}\right]^{-1}.
$$
Recall that this lower bound on $\gamma$ provides the upper bound on $\alpha_M = 1- \gamma$ which leads to the contradiction of the Lemma. Thus, by requiring $\alpha_M > 1 - \left[1+ \frac{c_1^n}{\epsilon^n}\right]^{-1}$, as originally prescribed in the assumption, we have proven the Lemma.
\end{proof}

{\bf Remark.}
Perelman's Maximal Volume Lemma proves the existence of a geodesic
in $M^n$ of length at least $c_1R > 1$ that is within a fixed
distance of a given point. Consider, for example, the case when $M^n
= \R^n$. Given a point $a \in \R^n$, it is possible to find a geodesic
of any length (in fact, there exists a ray) that is arbitrarily
close to $a$. Indeed, letting $c_1 \rightarrow \infty$ in the
expression for $\alpha_M$, while keeping $\epsilon$ and $n$ fixed,
we find that $\alpha_M \rightarrow 1$.  Similary, letting $\epsilon
\rightarrow 0$ (with $c_1, n$ fixed), forces $\alpha_M \rightarrow
1$ as well. Recall that by the Bishop-Gromov Volume Comparison Theorem,
$\alpha_M =1$ implies $M^n$ is isometric to $\R^n$.

{\bf Remark.}
Allowing the dimension of $M^n$ to increase while keeping $\epsilon$ and $c_1$
constant also pushes the lower bound on $\alpha_M$ closer to 1.

\section{Almost Equicontinuity and the Construction of Homotopies}\label{Chapter-almost equicontinuity}
In this section, we prove a general method of constructing
homotopies from sequences of increasingly refined nets. We begin by
reviewing a definition and theorem from \cite{So}.

\subsection{Background and Definitions}

\begin{Def}\label{Def-almost equicontinuous}
{\em [\cite{So}, Definition 2.5]} A sequence of functions between compact
metric spaces $f_i: X_i \to Y_i$, is said to be \textit{almost
equicontinuous} if there exists $\epsilon_i$ decreasing to 0 such
that for all $\epsilon > 0$ there exists $\delta_{\epsilon} >0$ such
that
\begin{equation}
d_{Y_i}(f_i(x_1), f_i(x_2)) < \epsilon + \epsilon_i,
~~~\textrm{whenever} ~~d_{X_i}(x_1, x_2) < \delta_{\epsilon}.
\end{equation}
\end{Def}

\begin{thm}\label{Theorem-Sormani}
{\em [\cite{So}, Theorem 2.3]} \textit{If $f_i:X_i \to Y_i$ is almost
equicontinuous between complete length spaces $(X_i, x_i) \to (X,x)$
and $(Y_i, y_i) \to (Y,y)$ which converge in the Gromov-Hausdorff
sense where $X$ and $Y$ are compact, then a subsequence of the $f_i$
converge to a continuous limit function $f: X \to Y$.}
\end{thm}

Let $X$ be a complete length space and let $K_j$ be a sequence of
finite cell decompositions of $X$. Each such decomposition $K_j$ is composed of a collection of cells $\sigma_i$ so that, for each $j$, $X = \coprod_{\sigma_i \in K_j} \sigma_i$. Each $K_{j+1}$ is a refinement of $K_j$.



\begin{Def}
Let $K$ be a finite cell decomposition of a complete length space
$X$. A map $\psi_K : X \to X$ which maps all the points in a cell $\sigma$ of $K$ to a single point $p \in \sigma$ is called a \textit{discrete
decomposition map of K}.
\end{Def}

\begin{Lem} \label{Lemma-Munn}
\textit{Let $K_j$ be a sequence of finite cell decompositions of $X$
and $\{ \psi_{K_j} \}$ a sequence of discrete decomposition maps of
$K_j$. This sequence of maps is almost equicontinuous provided
$\max\{\diam(\sigma) | \sigma \in K_j \} \to 0$ as $j \to \infty$.}
\end{Lem}

\begin{proof} For each $j$, let $d_j =\max\{\diam(\sigma) | \sigma \in K_j \}$.
Pick $\epsilon > 0$ and suppose $x,y \in X$ such that $d(x,y) <
\epsilon$. By the triangle inequality,  
\begin{eqnarray*}
d(\psi_j(x), \psi_j(y)) &\leq& d(\psi_j(x), x) + d(x, y) + d(\psi(y),y)\\
&<& \epsilon + 2d_j.
\end{eqnarray*}

Each $K_{j+1}$ is a refinement of $K_j$ and so by assumption the sequence $d_j$ decreases to 0. Thus, the sequence $\{ \psi_j \}$ is almost equicontinuous as claimed.
\end{proof}

\begin{Lem} \label{Lemma-Concat}
\textit{The composition of two almost equicontinuous sequences of
maps is again almost equicontinuous; i.e. if $\{f_j\}$ and $\{g_j\}$
are two sequences of maps which are almost equicontinuous. Then
$\{f_j \circ g_j\}$ is also almost equicontinuous.}
\end{Lem}

\begin{proof}
Suppose $\{f_j\}$ and $\{g_j\}$ are two almost equicontinuous
sequences of maps.   Since $\{f_j\}$ is almost equicontinous, given
$\epsilon>0$, there exists $\delta^f_{\epsilon}>0$ and positive
integer $K^f$ such that $d(f_j(x), f_j(y)) \leq \epsilon$ for all
$j> K^f$, provided $d(x,y) < \delta^f_{\epsilon}$. Choose $\delta^{f
\circ g}_{\epsilon} = \delta^f_{\delta^g_{\epsilon}}$ and choose a
positive integer $K = \max\{K^f, K^g\}$, where $\delta^g_{\epsilon}$
and $K^g$ are chosen so that when $d(a,b) < \delta^g_{\epsilon}$, we
have $d(g_j(a), g_j(b)) < \delta^f_{\epsilon}$, for all $j> K^g$.

Therefore, if $d(a,b) < \delta^f_{\delta^g_{\epsilon}}$, then
$d(g_j(a), g_j(b)) < \delta^f_{\epsilon}$, for all $j > K \geq K^g$
and thus, $d(f_j(g_j(a)), f_j(g_j(b))) < \epsilon$, for all $j
> K \geq K^f$. Therefore, the sequence $\{f_j \circ g_j\}$ is almost
equicontinuous.
\end{proof}

\subsection{Homotopy Construction Theorem}

The following theorem is crucial in constructing the homotopies in the manifold setting. In the statement of the theorem and in what follows we often refer to the $i$-skeleton of a cell decomposition $K$. We define an $i$-skeleton here.

\begin{Def}\label{Def-i-skeleton} The \textit{$i$-skeleton of a $k$-dimensional cell
decomposition K}, denoted $\skel_i(K)$ for $i=0,1,..,k$, is defined
as the collection of all $i$-dimensional cells contained in $K$.
\end{Def}

Note that if $X=\D^{k+1}$ then $\s^k \subset \D^{k+1}$ is contained
in $\skel_k(K)$ for any cell decomposition $K$ of $\D^{k+1}$.

\begin{thm}{\em \textbf{(Homotopy Construction Theorem).}} \label{Theorem-HomotopyConstruction}
\textit{Let $Y$ be a complete, locally compact metric space, $p \in
Y$, $R>0$ and $f: \s^k \to B_p(R) \subset Y$ a continuous map. Given
constants $c>1$, $\omega \in (0,1)$, and a sequence of finite cell
decompositions $K_j$ of $\D^{k+1}$ with maps $f_j: \skel_k(K_j) \to
Y$ satisfying
the following three properties\\
(A) $K_{j+1}$ is a subdivision of $K_j$ and $f_{j+1}
\equiv f_j$ on $K_j$ and $\max\{\diam(\sigma) | \sigma \in K_j\} \to 0$,\\
(B) For each $(k+1)$-cell, $\sigma \in K_j$, there exists a point
$p_{\sigma} \in B_p(cR) \subset Y$ and a constant $R_{\sigma} >0$
such that
$$f_j(\del \sigma) \subset B_{p_{\sigma}}(R_{\sigma});$$
and, if $\sigma' \subset \sigma$, where $\sigma'\in K_{j+1}$,
$\sigma \in K_j$, then $$B_{p_{\sigma'}}(cR_{\sigma'}) \subset
B_{p_{\sigma}}(cR_{\sigma}), \quad \textrm{and }~ R_{\sigma'} \leq
\omega R_{\sigma}, \textrm{for } \omega \in(0,1).$$ (C)
$\skel_k(K_0) = \s^k = \del \D^{k+1}$, $p_{\sigma_0}=p$, and
$R_{\sigma_0} = R$,\\ then the map $f$ can be continuously extended
to a map $g:\D^{k+1} \to B_p(cR) \subset Y$.}
\end{thm}

\begin{proof}
Suppose we have such a sequence of finite cell decompositions $K_j$
of $\D^{k+1}$ and continuous maps $f_j: \skel_k(K_j) \to M$
satisfying (A), (B), and (C) above. For any $x\in \D^{k+1}$, choose
a sequence of $(k+1)$-cells $\sigma_j \in K_j$, such that
$\sigma_{j+1} \subset \sigma_j$ and $x \in \clos(\sigma_j)$ for all
$j$. Therefore, each point $x \in \D^{k+1}$ determines a sequence of
$(k+1)$-cells `converging to' $x$. Each of these cells determines a
point, $p_{\sigma_j}$, and a radius, $R_{\sigma_j}>0$, which we
assume satisfy the properties outlined in (A), (B), and (C) above.

As in Perelman's homotopy construction \cite{P2}, define $g$ by
$g(x) = \lim_{j \to \infty} p_{\sigma_j}$. If $x\in \skel_k(K_j)$
for some $j$, set $g(x) = f_j(x)$. If $x \notin \skel_k(K_j)$ for all $j$, then for $j, k>0$, property (B) implies
\begin{equation*}
d(p_{\sigma_j},p_{\sigma_{j+k}}) \leq cR_{\sigma_j} \leq c\omega^j R.
\end{equation*}
Since $\omega \in (0,1)$, the sequence $\{p_{\sigma_j}\}$ is a
Cauchy sequence and thus converges. Hence, $g(x)$ is well-defined.

Note that $\del \D^{k+1} = \s^k = \skel_k(K_0)$ and so by the
definition of $g$, for any $x \in \del \D^{k+1}$,  $g(x) = f_0(x) =
f(x)$. Thus, $g|_{\del \D^{k+1}} = f$.

The continuity of $g$ is not verified in \cite{P2}. Here we prove
that $g$ is continuous. Define a sequence of maps $g_j: \D^{k+1} \to
Y$ by $g_j(x) = p_{\sigma_j}$ for each $j$.

\textbf{Claim.} The sequence of maps $\{g_j\}$ is uniformly almost
equicontinuous.

\textit{Proof of Claim.} Define a sequence of intermediate maps
$\psi_{K_j} : \D^{k+1} \to \D^{k+1}$, where $\psi_{K_j}$ is a
discrete decomposition map for $K_j$. Note that $\im(\psi_{K_j})$ is
a discrete metric space. Define $\overline{g_j} : \im(\psi_{K_j})
\to X$ in such a way that $\overline{g_j} = g_j
|_{\im(\psi_{K_j})}$.

By (A) we have that $\max\{\diam(\sigma) | \sigma \in K_j\} \to 0$
as $j \to \infty$. Therefore, the sequence of decomposition maps
$\psi_{K_j}$ is almost equicontinuous by Lemma \ref{Lemma-Munn}.

The maps $\overline{g_j}$ are discrete and thus the sequence
$\{\overline{g_j}\}$ is almost equicontinuous.

Since $g_j = \overline{g_j} \circ \psi_{K_j}$, by Lemma
\ref{Lemma-Concat}, the sequence of maps $\{g_j\}$ is also uniformly
almost equicontinuous. This completes the proof of the Claim.

Finally, by Theorem \ref{Theorem-Sormani} (see \cite{So} for proof),
the limiting map $g$ is continuous. This completes the proof of
Proposition \ref{Theorem-HomotopyConstruction}.
\end{proof}

\section{Double Induction Argument} \label{Section-General Argument}
In this section we use Perelman's double induction argument outlined
in section \ref{Chapter-Intro} to prove Theorem
\ref{Theorem-Main}. We introduce a collection of constants which are
defined inductively. We define them here as they are necessary for
the induction statements.

\begin{Def}\label{Def-C(k,n) constants}
For $k,n \in \N$ and $i=0,1,..,k$, define constants
$C_{k,n}(i)$ iteratively as follows:
\begin{equation}
C_{k,n}(i) = (16k)^{n-1}\left(1+10 C_{k,n}(i-1)^n + 3 + 10 C_{k,n}(i-1)\right), \quad i \geq 1
\end{equation}
and $C_{k,n}(0)=1$. We denote $C_{k,n} = C_{k,n}(k)$.
\end{Def}

\begin{Def}\label{Def-six* function}
Define a function
\begin{equation}
h_{k,n}(x) = \left[1-10^{k+2} C_{k,n}x\left(1+\frac{x}{2k} \right)^k\right]^{-1}.
\end{equation}
\end{Def}
This function $h_{k,n}$ has a vertical asymptote at $x=\delta_{k,n}$
for some small value $\delta_{k,n}>0$, where $10^{k+2}
\emph{C}_{k,n}\delta_{k,n}\left(1 + \frac{\delta_{k,n}}{2k}
\right)^k =1$. Note that $h_{k,n}: (0,\delta_{k,n}) \to (1,\infty)$
is a smooth, one-to-one, onto, increasing function. Thus
$h_{k,n}^{-1} : (1, \infty) \to (0, \delta_{k,n})$ is well-defined.

Toward proving Theorem \ref{Theorem-Main}, we need to build the
homotopy as described earlier. This requires control on the volume
growth of $M^n$. We now define the expression $\beta(k,c,n)$ which
we will use to control the volume growth of $M^n$.

\begin{Def}\label{Def-beta volume growth}
For constants, $c>1$ and $k,n \in \N$, the value of $\beta(k,c,n)$
represents a minimum volume growth necessary to guarantee that any
continuous map $f: \s^k \to B_p(R)$ has a continuous extension $g:
\D^{k+1} \to B_p(cR)$. Define
\begin{eqnarray}
\beta(k,c,n) &= \max\{&1-\gamma(c,h^{-1}_{k,n}(c), n);
\label{Def-AbGr for
beta}\\
&& \beta(j, 1 + \frac{h^{-1}_{k,n}(c)}{2k}, n), j=1,..,k-1
\label{Def-Main(k-1) for beta}\},
\end{eqnarray}
where $\beta(0,c,n)=0$ for any $c$ and $\beta(1,c,n) = 1-\gamma(c,
h^{-1}_{1,n}(c),n)$. Recall that $\gamma(c,d,n) =
[1+\frac{c^n}{d^n}]^{-1}$ [Definition \ref{Def-gamma}] was used in
proving Perelman's Maximal Volume Lemma [Lemma \ref{Lemma-Perelman's
Maximal Volume}].
\end{Def}

\subsection{Key Lemmas} \label{Section-Key Lemmas}
In this section we state the Main Lemma and the Moving In Lemma.
These are similar to the lemmas used in Perelman's paper \cite{P2}
except that we are controlling the constants carefully so as to
be able to determine the best bounds for the volume growth later.

\begin{Lem}  \emph{{\bf [Main Lemma($k$)].}} \label{Lemma-Main}
Let $M^n$ be a complete Riemannian manifold with $\Ric \geq 0$ and
let $p\in M^n$ and $R>0$. For any constant $c>1$ and $k,n \in \N$,
if
\begin{equation}\label{volume growth for Main(k)}
\alpha_M \geq \beta(k, c, n),
\end{equation}
then any continuous map $f: \s^k \to B_p(R)$ can be continuously
extended to a map $g: \D^{k+1} \to B_p(cR)$. \end{Lem}

\begin{Lem} \emph{\textbf{[Moving In Lemma($k$)].}} \label{Lemma-Moving In}
Let $M^n$ be a Riemannian manifold with $\Ric \geq 0$. For
any constant $d_0 \in (0, \delta_{k,n})$ and $k,n \in \N$ if
\begin{equation}\label{volume growth for Minor(k)}
\alpha_M \geq \beta(k, h_{k,n}(d_0), n),
\end{equation}
then given $q \in M^n$, $\rho >0$, a continuous map $\phi: \s^k
\rightarrow B_q(\rho)$ and a triangulation $T^k$ of $\s^k$ such that
$\diam(\phi(\Delta^k)) \leq d_0 \rho$ for all $\Delta^k \in T^k$,
there exists a continuous map $\tphi:\s^k \rightarrow
B_q((1-d_0)\rho)$ such that
\begin{equation}\label{diam conclusion from moving in}
\diam(\phi(\Delta^k) \cup \tphi(\Delta^k)) \leq 10^{-k-1} \left(1+
\frac{d_0}{2k}\right)^{-k} (1-h_{k,n}(d_0)^{-1}) \rho.
\end{equation}
\end{Lem}

In the next two sections we prove these lemmas. 

Before proceeding to the proofs, it is perhaps helpful to provide some insight to the main ideas behind the two lemmas above and how they are related to one another. The Moving In Lemma is, in some sense, the primary tool in constructing the homotopy. In fact, this lemma is precisely the point in the argument where the volume growth restriction is introduced. The new map $\tphi$ constructed in the Moving In Lemma is not necessarily homotopic the original map $\phi$; however, we require their images to be `close' in the manifold by controlling very carefully and uniformly the distance between the images of triangulations between the two maps. The proof is constructive and to construct a map with these properties requires large amount of volume growth in $M^n$. The Main Lemma provides a way of keeping track of the volume growth required to produce the homotopy. It's requirement on the volume growth arises only in the fact that it's proof requires an application of the Moving In Lemma in the same dimension.

The two lemmas are related to one another through the choice of the small constant $d_0$ in the Moving In Lemma, the constant $c > 1$ in the Main Lemma and the double induction argument relating the two. For example, taking $d_0$ very small in the Moving In Lemma weakens the restriction on the volume growth there. However, the Main Lemma is proven by induction using the constant $c = 1 + d_0/2k$ in lower dimensions. Taking $c$ very close to 1 in the Main Lemma ultimately forces the volume growth to be very large, close to 1. Contrarily, taking a much larger $d_0 < 1$ in the Moving In Lemma immediately forces the volume growth to be close to 1. The difficulty in determining optimal bounds (via this method) for the volume growth as stated in Theorem \ref{Theorem-Main} arises in finding the balance between these two competing lemmas and choosing the best constants $d_0$ and $c$.

In section \ref{Section-Proof of Main Lemma} we
prove Main Lemma($k$) assuming Moving In Lemma($k$) and Main
Lemma($j$) for $j=1,..,k-1$. In section \ref{Section-Proof of Moving In Lemma} we prove Moving In
Lemma($k$) assuming Main Lemma($i$), $i=0,..,k-1$. In section \ref{Section-Proof of Main Theorem}  we
apply these lemmas to prove Theorem \ref{Theorem-Main}. We begin by proving Main Lemma (0).

\begin{Lem} \emph{\textbf{[Main Lemma(0)].}} \label{Lemma-Main(0)}
Let $X$ be a complete length space and let $p\in X$, $R>0$. For any
constant $c>1$, any continuous map $f: \s^0 \to B_p(R) \subset X$
can be continuously extended to a map $g: D^1 \to B_p(cR) \subset
X$.
\end{Lem}

\begin{proof} The image $f(\s^0)$ consists of two points, $p_1, p_2 \in X$. Since $X$ is a
complete length space, it is possible to find length minimizing
geodesics $\sigma_i$ connecting $p_i$ to $p$, for $i=1,2$. Define
$g$ so that $\im(g) = \sigma_1 \cup \sigma_2$ and $g(-1) = p_1$ and
$g(1) = p_2$. Thus, $g$ is a continuous extension of the map $f$ and
by construction $\im(g) \subset B_p(cR) \subset X$.
\end{proof}

\subsection{Proof of Main Lemma($k$)}\label{Section-Proof of Main Lemma}

\begin{proof}
The proof is by induction on $k$. When $k=0$, the result follows
from Lemma \ref{Lemma-Main(0)}. No assumption on volume growth is
necessary. Assume now that Main Lemma($i$) holds for $i=1,..,k-1$:
Given any constants $c_i
>1$, a continuous map $f: \s^i \to B_p(R)$ has a continuous
extension to a map $g: \D^{i+1} \to B_p(c_iR)$ provided $\alpha_M
\geq \beta(i, c_i, n)$. We will now show that the result is true for
dimension $k$.

Let $f: \s^k \to B_p(R)\subset M^n$ be a continuous map. Choose
$c>1$ and suppose $\alpha_M \geq \beta(k,c,n)$. Our goal now is to
show that the map $f: \s^k \to B_p(R)$ has a continuous extension.
To do this we will show that there exists a sequence of finite cell
decompositions, $K_j$, of $\D^{k+1}$ and maps $f_j$ that satisfy the
hypothesis of the Homotopy Construction Theorem [Theorem
\ref{Theorem-HomotopyConstruction}] and thus create the homotopy $g:
\D^{k+1} \to B_p(cR)$.

For $j=0$, define $K_0$ to be the cell decomposition consisting of a
single cell (i.e. $K_0 \cong \D^{k+1}$) so $\skel_k(K_0) = \s^k$.
Recall that we use the notation $\skel_k(K_j)$ to denote the union
of the boundaries of the cell decomposition of $K_j$ [Definition
\ref{Def-i-skeleton}].

As in \cite{P2}, inductively define $K_{j+1}$ given $K_j$ in the
following way. For a $(k+1)$-cell, $\sigma \in K_j$, note that
$\sigma$ is homeomorphic to a disk so it can be viewed in polar
coordinates as $\left(\s^k \times (0,1]\right) \cup \{0\}$. Let
$T^k_{\sigma}$ be a triangulation of $\s^k$, where $\s^k \cong
\del\sigma$ and $\diam_{\sigma}(\Del^k) < 1/k$ for all $\Del^k \in
T^k_{\sigma}$. Define $K_{j+1}$ so that
\begin{equation}
\sigma \cap \skel_k(K_{j+1}) = (\s^k \times \{1\}) \cup (\s^k \times
\{1/2\}) \cup \left(\skel_{k-1}(T^k_{\sigma}) \times[1/2,1]\right).
\end{equation}
This inductive construction of the $K_j$ provides us with a sequence
of finite cell decompositions of $\D^{k+1}$.  Note that with an
appropriate selection of $\s^k \times \{1/2\}$ this sequence of
decompositions satisfies Condition A on cell decompositions as
required by the Homotopy Construction Theorem [Theorem
\ref{Theorem-HomotopyConstruction}] because $\max \{\diam(\sigma) |
\sigma \in K_j\} \rightarrow 0$.

Next, we define the continuous maps $f_j: \skel_k(K_j) \to M^n$.
Begin by setting $f_0 \equiv f$. In this way, $f_0: \skel_k(K_0) \to
B_p(R) \subset M^n$ and the initializing hypothesis (C) of Theorem
\ref{Theorem-HomotopyConstruction} is satisfied. We verify the rest
of the hypothesis inductively.

Suppose $f_j$ satisfies hypotheses (A) and (B) of Theorem
\ref{Theorem-HomotopyConstruction}. It remains to define $f_{j+1}$
and check that hypotheses (A) and (B) hold for this $f_{j+1}$. We
describe the process to define $f_{j+1}$ on the refinement of a
single $(k+1)$-cell $\sigma \in K_j$. To define $f_{j+1}$ on all of
$\skel_k(K_{j+1})$, repeat this process on each $(k+1)$-cell of
$K_j$.

Given a $(k+1)$-cell $\sigma \in K_j$, by hypothesis (B), there
exists a point $p_{\sigma} \in B_p(cR) \subset M^n$ and a constant
$R_{\sigma}>0$ such that $f_j(\del \sigma) \subset
B_{p_{\sigma}}(R_{\sigma})$. As before, view $\sigma$ as $\left(\s^k
\times (0,1]\right) \cup \{0\}$, and think of $f_j$ as a map $f_j :
\s^k \to B_{p_{\sigma}}(R_{\sigma})$.

Define $f_{j+1} : \skel_k(K_{j+1}) \to M^n$ in three stages.

First we set
\begin{equation} \label{fedge1}
f_{j+1} \equiv f_j \qquad \textrm{on } \s^k \times \{1\},
\end{equation}
which is all that is required to satisfy hypothesis (A).

We claim that we can apply the Moving In Lemma($k$) to the map
$f_j$. Set $d_0 = h^{-1}_{k,n}(c)$ and keep $k,n$ as before. The
volume growth assumption (\ref{volume growth for Minor(k)}) is
satisfied since $\alpha_M \geq \beta(k,c,n) = \beta(k, h_{k,n}(d_0),
n)$.

Take $q=p_\sigma$, $\rho=R_\sigma$, and $\phi=f_j$ and take a
sufficiently fine triangulation, $T^k_{\sigma}$, of $\s^k \cong \del
\sigma$ such that $\diam(f_j(\Del^k)) \leq d_0 R_{\sigma}$ for all
$\Del^k \in T^k_{\sigma}$. Applying the Moving In Lemma($k$) [Lemma
\ref{Lemma-Moving In}], we obtain a map $\tilde{f_j} : \s^k \to
B_{p_{\sigma}}((1-d_0)R_{\sigma})$. We set
\begin{equation} \label{fedge2}
f_{j+1} \equiv \tilde{f_j} \qquad \textrm{on } \s^k \times \{1/2\}.
\end{equation}
This completes the second stage of our construction of $f_{j+1}$.
Furthermore, by (\ref{diam conclusion from moving in}),
\begin{eqnarray}
\hskip -.4in \diam(f_j(\Del^k) \cup \tilde{f_j}(\Del^k)) &\leq&
10^{-k-1}\left(1+\frac{d_0}{2k}\right)^{-k}(1-(h_{k,n}(d_0))^{-1})R_{\sigma}\\
&=&\label{usingmoving}
10^{-k-1}\left(1+\frac{d_0}{2k}\right)^{-k}(1-c^{-1})R_{\sigma},
\end{eqnarray}
for all $\Del^k \in T^k_{\sigma}$.

For the third stage and to complete the definition of $f_{j+1}$ on
$\sigma \cap \skel_k(K_{j+1})$, it remains to define $f_{j+1}$ on
$\skel_i(T^k_{\sigma}) \times [1/2, 1]$ for $i=0,1,..,k-1$. Below we
describe this procedure (inductively) for a single $k$-simplex
$\Del^k$ of the triangulation $T^k_{\sigma}$. Here we use the
induction hypothesis and assume the Main Lemma($j$) is true for
$j=1,..,k-1$. First, we apply Lemma \ref{Lemma-Main(0)} to the
0-skeleton [note that Lemma \ref{Lemma-Main(0)} is an analog of Main
Lemma(0)]. Then, we apply Main Lemma [Lemma \ref{Lemma-Main}]
repeatedly starting with the 1-dimension skeleton and continuing to
the $(k-1)$-dimension skeleton.

Let $\Del^0 \in T^k_{\sigma}$ be a 0-simplex. Consider the map
$f_{j+1,0}$ on $\s^0$ defined by $f_{j+1,0}(-1)= f_{j+1}(\Del^0
\times \{1\})$ and $f_{j+1,0}(1)= f_{j+1}(\Del^0 \times \{1/2\})$.
On these components, the map $f_{j+1,0}$ is obtained from
(\ref{fedge1}) and (\ref{fedge2}). We want to define $f_{j+1}$ on
$\Del^0 \times [1/2, 1]$. Note that,
\begin{eqnarray}
\hskip -.2in \diam(\im(f_{j+1,0})) &=& d(f_{j+1,0}(-1), f_{j+1,0}(1))\\
&=& \diam(f_{j+1}(\Del^0 \times \{1\}) \cup f_{j+1}(\Del^0 \times \{1/2\}))\\
&=& \diam(f_j(\Del^0) \cup \tilde{f_j}(\Del^0))\\
&\leq& \diam(f_j(\Del^k) \cup \tilde{f_j}(\Del^k))\\
&\leq&
10^{-k-1}\left(1+\frac{d_0}{2k}\right)^{-k}(1-c^{-1})R_{\sigma}.
\end{eqnarray}
In this last line we have applied (\ref{usingmoving}).

If we set
\begin{equation}
R_{j+1,0}=1/2 \cdot
10^{-k-1}\left(1+\frac{d_0}{2k}\right)^{-k}(1-c^{-1})R_{\sigma},
\end{equation}
then, by our estimate on the diameter of its image, we have
\begin{equation}
f_{j+1,0}: \s^0 \to B_{p_{j+1,0}}(R_{j+1,0}),
\end{equation}
for some point $p_{j+1,0}\in M^n$. We now apply Main Lemma(0) [Lemma
\ref{Lemma-Main(0)}] taking $c = 1+ d_0/2k$, $p=p_{j+1,0}$,
$R=R_{j+1,0}$ and $f=f_{j+1,0}$. Clearly, the hypotheses of Main
Lemma(0) are satisfied since $\beta(0,c,n)=0$ and $M^n$ is a
complete Riemannian manifold. Therefore, there exists a continuous
extension
\begin{equation}
g_{j+1,1}: \D^{1} \to
B_{p_{j+1,0}}\left(\left(1+\frac{d_0}{2k}\right)R_{j+1,0}\right)
\end{equation}
and we use it to define $f_{j+1}$ on $\skel_0(T^k_\sigma) \times
[1/2,1]$. Furthermore,
\begin{eqnarray}
\hskip -.5in \diam(f_{j+1}(\Del^0 \times [1/2,
1])) &=& \diam(\im(g_{j+1,1}))\\
&\leq& 2 \cdot \left(1+\frac{d_0}{2k}\right)R_{j+1,0}\\
&\leq& 2 \cdot
\left(1+\frac{d_0}{2k}\right)\cdot 1/2 \cdot\\
&& \qquad \left(10^{-k-1}\left(1+\frac{d_0}{2k}\right)^{-k}(1-c^{-1})R_{\sigma}\right)\\
&\leq&
10^{-k-1}\left(1+\frac{d_0}{2k}\right)^{-k+1}(1-c^{-1})R_{\sigma}.
\label{diam bound Del0}
\end{eqnarray}

We will use induction on $i$ to define $f_{j+1}$ on $\Del^i \times
[1/2, 1]$, for $0 \leq i < k$. Assume we have defined $f_{j+1} =
f_j$ on all simplices $\Del^i \in T^k_{\sigma}$ and we have defined
$f_{j+1}$ on all possible $\Del^{i-1} \times [1/2,1]$ so that
\begin{equation}
\diam(f_{j+1}(\Del^{i-1} \times [1/2, 1])) \leq
10^{i-1-k}\left(1+\frac{d_0}{2k}\right)^{i-k}(1-c^{-1})R_{\sigma}.\label{f_j+1
induction}
\end{equation}
Note that this holds for $i=1$ by (\ref{diam bound Del0}). Also,
note that (\ref{usingmoving}) implies
\begin{eqnarray}
\diam(f_{j+1}(\Del^i \times \{1\}) \cup f_{j+1}&& \hskip -.4in (\Del^i \times
\{1/2\}))\\
\hskip 1in &=& \diam(f_{j}(\Del^i) \cup \tilde{f}_{j}(\Del^i))\\
&\leq& \diam(f_{j}(\Del^k) \cup \tilde{f}_{j}(\Del^k))\\
&\leq&
10^{-k-1}\left(1+\frac{d_0}{2k}\right)^{-k}(1-c^{-1})R_{\sigma}.
\label{f_j+1 using moving}
\end{eqnarray}

We now build a new map $f_{j+1,i+1}$ on $\Del^i \times [1/2,1]$.
View\\ $(\Del^i \times \{1\}) \cup (\Del^i \times \{1/2\}) \cup
(\del\Del^i \times [1/2,1])$ as $\s^i$. Since $\del \Del^i \times
[1/2,1]$ is a collection of $\Del^{i-1} \times [1/2,1]$, we have a
map
\begin{equation}
f_{j+1,i}: \s^i \to B_{p_{j+1,i}}(R_{j+1,i}),
\end{equation}
for some point $p_{j+1,i}\in M^n$ and where by (\ref{f_j+1
induction}) and (\ref{f_j+1 using moving}) we have

\begin{eqnarray}
2R_{j+1,i} &=& \diam(f_{j+1}|_{\Del^i \times \{1\}} \cup
f_{j+1}|_{\Del^i \times \{1/2\}}) + \\
&& \hskip 1.5in
\diam(f_{j+1}(\del\Del^{i}\times [1/2,1]))\\
&\leq& \diam(f_{j}(\Del^i) \cup \tilde{f_{j}}(\Del^i)) +\\
&& \hskip 1.5in
\diam(f_{j+1}(\Del^{i-1}\times [1/2,1]))\\
&\leq&
10^{-k-1}\left(1+\frac{d_0}{2k}\right)^{-k}(1-c^{-1})R_{\sigma} + \\
&& \hskip 1.25in 10^{i-1-k}\left(1+\frac{d_0}{2k}\right)^{i-k}(1-c^{-1})R_{\sigma}\\
&\leq&
10^{i-k}\left(1+\frac{d_0}{2k}\right)^{-k+i}(1-c^{-1})R_{\sigma}.
\end{eqnarray}

Therefore,
\begin{equation}
\diam(\im(f_{j+1,i})) \leq
10^{i-k}\left(1+\frac{d_0}{2k}\right)^{-k+i}(1-c^{-1})R_{\sigma}.
\end{equation}
Apply Main Lemma($i$) taking $c=1+d_0/2k$ and $k,n$ as before. This
is allowed because the volume growth requirement for Main Lemma($i$)
is satisfied by (\ref{Def-Main(k-1) for beta}) and because the volume growth satifies
\begin{eqnarray}
\alpha_M &\geq& \beta(k,c,n)\\
&\geq& \beta(i, 1+\frac{h^{-1}_{k,n}(c)}{2k},n)\\
&=& \beta(i, 1+\frac{d_0}{2k},n).
\end{eqnarray}
Therefore, there exists a continuous extension
\begin{equation}
g_{j+1,i+1} : \D^{i+1} \to B_{p_{j+1,i}}((1+d_0/2k)R_{j+1,i})
\end{equation}
of the continuous map $f_{j+1,i}$. This extension defines $f_{j+1}$
on $\skel_{i}(T^k_{\sigma}) \times [1/2,1]$ and we have the bound
\begin{eqnarray}
\hskip -.4in \diam(f_{j+1}(\Del^i \times [1/2,1])) &=&
\diam(\im(g_{j+1,i+1}))\\
&\leq& 2 \cdot \left(1+\frac{d_0}{2k}\right) \cdot R_{j+1,i}\\
&=& 2 \cdot \left(1+\frac{d_0}{2k}\right) \cdot 1/2 \cdot\\
&& \hskip.5in 10^{i-k}\left(1+\frac{d_0}{2k}\right)^{-k+i}(1-c^{-1})R_{\sigma}\\
&=&
10^{i-k}\left(1+\frac{d_0}{2k}\right)^{-k+i+1}(1-c^{-1})R_{\sigma}.
\end{eqnarray}
Furthermore, we have the bound
\begin{equation} \label{diam bound for i simplex}
\diam (f_{j+1}(\Del^i \times [1/2,1])) \leq
10^{i-k}\left(1+\frac{d_0}{2k}\right)^{i+1-k}(1-c^{-1})R_{\sigma},
\end{equation}
for all $\Del^i \subset \Del^k$, $i=0,1,..,k-1$, which implies our
induction hypothesis on $i$. Thus, we have defined $f_{j+1}$ on
$\skel_i(T^k_{\sigma}) \times [1/2,1]$ for each $i=0,1,..k-1$.

We now complete the proof by showing that the hypotheses (A) and (B)
of the Homotopy Construction Theorem [Theorem
\ref{Theorem-HomotopyConstruction}] hold for the function $f_{j+1}$.

Hypothesis (A) follows immediately from this construction since each
$K_{j+1}$ is a subdivision of the previous $K_j$ and by definition
$f_{j+1} \equiv f_j$ on $K_j$.

To check (B) holds, let $\sigma' \in K_{j+1}$ and suppose $\sigma'
\cong \Del^k \times [1/2,1]$ for some $\Del^k \in \s^k$. Notice that
\begin{eqnarray}
\hskip -.4in \diam(f_{j+1}(\del \sigma')) &\leq& \diam(f_{j+1}|_{\Del^k \times
\{1\}} \cup f_{j+1}|_{\Del^k \times \{1/2\}}) +\\
&& \hskip 1.5in \diam(f_{j+1}(\del\Del^{k}\times [1/2,1]))\\
&\leq& \diam(f_{j}(\Del^k) \cup \tilde{f_{j}}(\Del^k)) +\\
&& \hskip 1.4in \diam(f_{j+1}(\Del^{k-1}\times [1/2,1]))\\
&\leq&
10^{-k-1}\left(1+\frac{d_0}{2k}\right)^{-k}(1-c^{-1})R_{\sigma} +\\
&& \hskip 2in 10^{-1}(1-c^{-1})R_{\sigma},\label{diam bound on f_j+1}
\end{eqnarray}
where the last line follows from (\ref{usingmoving}) and (\ref{diam
bound for i simplex}) with $i=k-1$.

Set
\begin{equation}\label{omega}
R_{\sigma'} = 1/2 \cdot
[10^{-k-1}\left(1+\frac{d_0}{2k}\right)^{-k}(1-c^{-1}) +
10^{-1}(1-c^{-1})]R_{\sigma}.
\end{equation}
Then, by (\ref{diam bound on f_j+1}), there exists a point
$p_{\sigma'} \in M^n$ such that $f_{j+1}(\del\sigma') \subset
B_{p_{\sigma'}}(R_{\sigma'}).$

To verify $B_{p_{\sigma'}}(cR_{\sigma'}) \subset
B_{p_{\sigma}}(cR_{\sigma})$, let $x \in
B_{p_{\sigma'}}(cR_{\sigma'})$ and notice that for $q \in f(\Del^k
\times \{1/2\}) \subset B_{p_{\sigma'}}(R_{\sigma'})$,
\begin{eqnarray}
d(x, p_{\sigma}) &\leq& d(x,q) + d(q,
p_{\sigma})\\
&\leq& 2 \cdot 1/2 (1-c^{-1}) c R_{\sigma} +
(1-d_0)R_{\sigma}\\
&\leq& (c-1) R_{\sigma} + (1-d_0)R_{\sigma}\\
&<& cR_{\sigma}.
\end{eqnarray}
Therefore, $B_{p_{\sigma'}}(cR_{\sigma'}) \subset
B_{p_{\sigma}}(cR_{\sigma})$.

Furthermore, since $B_{p_{{\sigma'}}}(cR_{\sigma'}) \subset
B_{p_{\sigma}}(cR_{\sigma})$ for all nested sequences $\sigma'
\subset \sigma$, it follows that
\begin{eqnarray}
d(p_{\sigma'},p) &\leq& d(p_{\sigma'}, p_{\sigma}) + ... + d(p_{\sigma_.},p)\\
&\leq& cR_{\sigma} - cR_{\sigma'} +...+ cR - cR_{\sigma_.}\\
&=& cR - cR_{\sigma'}\\
&<& cR.
\end{eqnarray}
Thus, $p_{\sigma'} \in B_p(cR)$ as required.

Lastly, we have $R_{\sigma'} \leq \omega R_{\sigma}$ for
\begin{equation}
\omega = 1/2 \cdot \left[10^{-k-1}
\left(1+\frac{d_0}{2k}\right)^{-k}(1-c^{-1}) + 10^{-1}(1-c^{-1})\right].
\end{equation}
Note that $\omega \in (0,1)$ because $k \geq 1$ and $d_0<1$.

Thus, we have constructed a sequence of maps $f_j: \skel_k(K_j) \to
M^n$ satisfying the hypotheses of the Homotopy Construction Theorem
[Theorem \ref{Theorem-HomotopyConstruction}].
Therefore, the map $f$ can be continuously extended to a map\\ $g:
\D^{k+1} \to B_p(cR) \subset M^n$. This completes the proof of Main
Lemma($k$).
\end{proof}

\subsection{Proof of Moving In Lemma($k$)}\label{Section-Proof of Moving In Lemma}
We now prove Moving In Lemma($k$) assuming that Main Lemma($j$) is
true for $j=0,..,k-1$.

\begin{proof}Recall that $\alpha_M \geq \beta(k, h_{k,n}(d_0),n)$ and
we are given $q \in M^n$, $\rho >0$, a continuous map $\phi: \s^k
\rightarrow B_q(\rho)$ and a triangulation $T^k$ of $\s^k$ such that
$\diam(\phi(\Delta^k)) \leq d_0 \rho$ for all $\Delta^k \in T^k$. We
must show that there exists a continuous map $\tphi:\s^k \rightarrow
B_q((1-d_0)\rho)$ such that
\begin{equation}\label{diam bound from moving in}
\diam(\phi(\Delta^k) \cup \tphi(\Delta^k)) \leq 10^{-k-1} \left(1+
\frac{d_0}{2k}\right)^{-k} (1-h_{k,n}(d_0)^{-1}) \rho.
\end{equation}

We will construct $\tphi$ inductively on $\skel_i(T^k)$ for $i =
0,..,k$ in such a way that $\tphi(\Del^i)) \equiv \phi(\Del^i)$ if
$\phi(\Del^i) \subset B_q((1-2d_0)\rho)$; and, if $\phi \nsubseteq
B_q((1-2d_0)\rho)$, then
\begin{eqnarray}
\tphi(\Del^i) &\subset& B_q((1-d_0(2-i/k))\rho),\label{cond1} \\
\diam(\phi(\Del^i) \cup \tphi(\Del^i)) &\leq& 10 d_i \rho,
\label{cond2}
\end{eqnarray}
for all $\Del^i \subset T^k$, $i=0,..,k$. The constants $d_i>0$
satisfy

\vskip -.25in

\begin{eqnarray}
d_0 + 10d_i &\leq& b_i(d_{i+1} - 3d_0 - 10d_i)\label{Ineq1}\\
d_0 + 10d_i &\leq& b_i(c-1 + d_0(2-i/k))\label{Ineq2}\\
8 b_i^{\frac{1}{n-1}}(d_0 + 10d_i) &\leq& \frac{d_0}{2k}\label{Ineq3}\\
10d_k &\leq&
10^{-k-1}(1+d_0/2k)^{-k}(1-h_{k,n}(d_0)^{-1}),\label{Ineq4}
\end{eqnarray}
for some constants $b_i \in (0,1/2]$. The existence of such
constants $d_i$ and $b_i$ is proven in Lemma \ref{Lemma-d_i and
b_i}. Note that (\ref{cond2}) and (\ref{Ineq4}) together immediately imply
(\ref{diam bound from moving in}). Thus, we need only define $\tphi$ so that
the above conditions are obeyed. To do so, we construct $\tphi$ successively
on the $i$-skeleta of $T_k$.

Begin with the case $i=0$. Let $\Del^0 \in \skel_0(T^k)$ and assume
$\phi(\Del^0) \notin B_q((1-2d_0)\rho)$, else we are done. Let
$\sigma_{\Del^0}$ denote a length minimizing geodesic from
$\phi(\Del^0)$ to $q$ and define $\tphi(\Del^0) =
\sigma_{\Del^0}((1-2d_0)\rho)$. In this way, $\tphi(\Del^0) \in
B_q((1-2d_0)\rho)$ and (\ref{cond1}) is satisfied for $i=0$.
Furthermore,
\begin{eqnarray}
\diam(\phi(\Del^0) \cup \tphi(\Del^0)) &=& d(\phi(\Del^0),
\tphi(\Del^0))\\
&=&d(q,\phi(\Del^0)) - d(q, \tphi(\Del^0))\\
&\leq& \rho - (1-2d_0)\rho = 2d_0 \rho \leq 10 d_0 \rho.
\end{eqnarray}
Thus, (\ref{cond2}) is also satisfied when $i=0$.

Now assume that $\tphi$ is defined on $\skel_i(T^k)$ and that
(\ref{cond1}) and (\ref{cond2}) for $0\leq i
\leq k-1$. We now construct $\tphi$ on $\skel_{i+1}(T^k)$. Let $\Del^{i+1}
\subset \skel_{i+1}(T^k).$ As before, suppose $\phi(\Del^{i+1})
\nsubseteq B_q((1-2d_0)\rho)$, else we are done by simply setting
$\tphi(\Del^{i+1}) \equiv \phi(\Del^{i+1})$.

Next apply Perelman's Maximal Volume Lemma [Lemma
\ref{Lemma-Perelman's Maximal Volume}], taking $c_1=h_{k,n}(d_0)$,
$\epsilon = d_0$, and $p=q$, $R=\rho$. Since, by our hypothesis,
\begin{eqnarray}
\hskip -.8in \alpha_M \hskip -.1in &\geq& \beta(k, h_{k,n}(d_0),n)\\
&=& \hskip -.1in \max \left\{1-\gamma(h_{k,n}(d_0),d_0, n); \beta\left(j, 1 +
\frac{d_0}{2k}, n\right), j=1,..,k-1 \right\}\\
&\geq& 1-\gamma(h_{k,n}(d_0),d_0, n),
\end{eqnarray}
there exists a point $r_{\Del} \in M^n \setminus B_q(h_{k,n}(d_0)
\rho)$ such that $d(\phi(\Del^{i+1}),\overline{q r_{\Del}}) \leq d_0
\rho$. Recall, $\overline{q r_{\Del}}$ denotes a minimal geodesic
connecting $q$ and $r_{\Del}$. Let $\sigma_{\Del}$ be a length
minimizing geodesic from $q$ to $r_{\Del}$ and define a point
$q_{\Del} = \sigma_{\Del}((1-d_{i+1})\rho)$. For any $x\in
\del\Del^{i+1}$, the triangle with vertices $\tphi(x)$, $q_{\Del}$,
and $r_{\Del}$ is small and thin. To verify this, we use the the
induction hypothesis that $\tphi$ has already been defined on
$\skel_i(T^k)$, $0\leq i \leq k-1$, and that the properties
(\ref{cond1}),(\ref{cond2}) are satisfied in dimension $i$.

Note that,
\begin{eqnarray}
d(\tphi(x), \overline{q_{\Del}r_{\Del}}) &\leq& d(\phi(x), \overline{q_{\Del}}r_{\Del}) + d(\phi(x),\tphi(x))\\
&\leq& d_0 \rho + \diam(\phi(\Del^i) \cup \tphi(\Del^i))\\
&\leq& d_0 \rho + 10 d_i \rho.
\end{eqnarray}
And
\begin{eqnarray}
\hskip -.5in d(\tphi(x), q_{\Del}) &\geq& d(q_{\Del},\phi(x)) -d(\phi(x),\tphi(x))\\
&\geq& d(q,\phi(x)) - d(q_{\Del},q) - \diam(\phi(\Del^i) \cup \tphi(\Del^i))\\
&\geq& (1-2d_0)\rho - \diam(\phi(\Del^{i+1}) - (1-d_{i+1})\rho -10d_i\rho\\
&\geq& (1-2d_0) \rho -d_0 \rho -(1-d_{i+1})\rho - 10d_i \rho\\
&\geq& (d_{i+1} - 3d_0 - 10d_i)\rho.
\end{eqnarray}
And finally,
\begin{eqnarray}
d(\tphi(x), r_{\Del}) &\geq& d(r_{\Del}, q) - d(q, \tphi(x))\\
&\geq& C\rho - d(q, \tphi(\Del^i))\\
&\geq& c\rho - (1-d_0(2-i/k))\rho\\
&=& (c-1 + d_0(2-i/k))\rho.
\end{eqnarray}
The inequalties (\ref{Ineq1}) and (\ref{Ineq2}) guarantee that the
triangle with vertices $\tphi(x)$, $q_{\Del}$, and $r_{\Del}$ is
small and thin for some constants $0 < b_i \leq 1/2$.

According to the excess estimate of Abresch-Gromoll [Theorem
\ref{Theorem-Abresch-Gromoll Excess}], we have that for any $x \in
\del \Del^{i+1}$, with $i=0,1,..,k-1$,
\begin{eqnarray}
\hskip -.4in e_{q_{\Del}, r_{\Del}}(\tphi(x)) &=& d(\tphi(x), q_{\Del}) + d(\tphi(x), r_{\Del}) - d(q_{\Del},r_{\Del})\\
&\leq& 8 \left(\frac{d(\tphi(x),
\overline{q_{\Del}r_{\Del}})}{\min\{d(\tphi(x), q_{\Del}),
d(\tphi(x), r_{\Del})\}}\right)^{\frac{1}{n-1}}
d(\tphi(x), \overline{q_{\Del}r_{\Del}})\\
&\leq& 8 b_{i}^{\frac{1}{n-1}}(d_0 + 10d_i)\rho.\label{excess}
\end{eqnarray}

Also, by the triangle inequality,
\begin{eqnarray}\label{triangle neq}
d(q,q_{\Del}) + d(q_{\Del},r_{\Del}) = d(q, r_{\Del}) \leq d(q,
\tphi(x)) + d(\tphi(x), r_{\Del}).
\end{eqnarray}

Adding (\ref{excess}) and (\ref{triangle neq}), we get
\begin{eqnarray}
\hskip -.4in d(\tphi(x), q_{\Del}) &\leq& 8 b_i^{\frac{1}{n-1}}(d_0 + 10d_i)\rho
+ d(q_{\Del}, r_{\Del}) - d(\tphi(x), r_{\Del})\\
&\leq& 8 b_i^{\frac{1}{n-1}}(d_0 + 10d_i)\rho + d(q,
\tphi(x)) + d(\tphi(x), r_{\Del})\\
&& \hskip 1.8in - d(q,q_{\Del}) - d(\tphi(x), r_{\Del})\\
&\leq& 8 b_i^{\frac{1}{n-1}}(d_0 + 10d_i)\rho + (1-d_0(2-i/k))\rho -
(1-d_{i+1}\rho)\\
&=& \left(8 b_i^{\frac{1}{n-1}}(d_0 + 10d_i) + d_{i+1} -
d_0(2-i/k))\right)\rho
\end{eqnarray}

It then follows from (\ref{Ineq3}) that, for all $x \in
\del\Del^{i+1}$,
\begin{equation}
d(\tphi(x), q_{\Del}) \leq (\frac{d_0}{2k} + d_{i+1} d_0(2-i/k))\rho = \left(d_{i+1} -d_0(2-\frac{2i+1}{2k}) \right)\rho.
\end{equation}
Now apply the Main Lemma [Lemma \ref{Lemma-Main}] in dimension $i$
taking
\begin{eqnarray}
p&=&q_{\Del},\\
R&=&\left(d_{i+1} - d_0 \left(2 - \frac{2i+1}{2k} \right) \right)\rho,\\
c&=&1+d_0/2k;
\end{eqnarray}
and letting $f=\tphi$. Since
\begin{eqnarray}
\hskip -.8in \alpha_M &\hskip -.1in \geq& \beta(k, h_{k,n}(d_0),n)\\
&\hskip -.1in =& \hskip -.1in \max \hskip -.05in \left\{1-\gamma(h_{k,n}(d_0),d_0, n); \beta\left(j, 1 +
\frac{d_0}{2k}, n\right), j=1,..,k-1 \right\}\\
&\hskip -.1in \geq& \beta(i, 1 + \frac{d_0}{2k}, n)
\end{eqnarray}
by our hypothesis, there exists a continuous extension of $\tphi$
from $\del \Del^{i+1}$ to $\Del^{i+1}$. Furthermore,
\begin{eqnarray}
d(\tphi(\Del^{i+1}), q_{\Del}) &\leq& (1+d_0/2k)\left( d_{i+1}
-d_0(2-\frac{2i+1}{2k})\right)\\
&\leq& \left(d_{i+1} -d_0(2- \frac{i+1}{k}) \right)\rho,
\end{eqnarray}
provided $d_i <1$, which is guaranteed by the fact that the $d_i$'s
are increasing in $i$ and, by (\ref{Ineq4}), $d_k <1$. Therefore, by the
triangle inequality,
\begin{eqnarray}
d(\tphi(\Del^{i+1}), q) &\leq& d(\tphi(\Del^{i+1}), q_{\Del}) +
d(q_{\Del}, q)\\
&\leq& \left(d_{i+1} -d_0(2- \frac{i+1}{k})\right)\rho + \left(1 -
d_{i+1}\right)\rho\\
&=& \left(1 - d_0(2- \frac{i+1}{k})\right)\rho.
\end{eqnarray}

Thus, (\ref{cond1})  is satisfied for $i+1$ for any choice of $d_i$,
$b_i$ satisfying the inequalities (\ref{Ineq1}), (\ref{Ineq2}) and
(\ref{Ineq3}).

Furthermore,
\begin{eqnarray}
\hskip -.5in \diam(\phi(\Del^{i+1}\cup\tphi(\Del^{i+1})) &\leq& \diam(\phi(\del\Del^{i+1}) \cup \tphi(\del\Del^{i+1})) +\\
&& \hskip0.5in \diam(\phi(\Del^{i+1})) + \diam(\tphi(\Del^{i+1}))\\
&\leq& \hskip -.05in 10d_i \rho + d_0\rho + 2\left(d_{i+1} -d_0(2- \frac{i+1}{k}\right)\rho\\
&=& \hskip -.05in \left( 2d_{i+1} + d_0\left(-3 + \frac{2(i+1)}{k}\right) +
10d_i\right) \rho\\
&\leq& (2d_{i+1} + 10d_i - d_0) \rho.
\end{eqnarray}
The inequality (\ref{Ineq1}) and the fact that $0 < b_i \leq 1/2$ imply that
\begin{equation}
\diam(\phi(\Del^{i+1})\cup\tphi(\Del^{i+1})) \leq 10d_{i+1}\rho.
\end{equation}

So, (\ref{cond2}) are satisfied for dimension $i+1$. Therefore,
$\tphi$ has been defined so that (\ref{cond1}) and (\ref{cond2}) are
satisfied for $i=0,..,k$. When $i=k$, (\ref{cond1}) implies
\begin{equation}
\tphi(\Del^k) \subset B_q((1-d_0)\rho), \qquad \forall \Del^k \in
T^k.
\end{equation}
Thus, we have constructed the map $\tphi: \s^k \to
B_q((1-d_0)\rho)$; and furthermore,
\begin{eqnarray}
\hskip -.4 in\diam(\phi(\Del^k) \cup \tphi(\Del^k)) &\leq& 10d_k \rho\\
&\leq& 10^{-k-1}\left(1 + \frac{d_0}{2k} \right)^{-k}
\left(1-h_{k,n}(d_0)^{-1} \right),
\end{eqnarray}
where the last inequality follows from (\ref{cond2}).
\end{proof}

\subsection{Proof of the Main Theorem} \label{Section-Proof of Main Theorem} In this section we
prove Theorem \ref{Theorem-Main} using Main Lemma($k$).

As a direct consequence of Main Lemma($k$) [Lemma \ref{Lemma-Main}],
we have

\begin{Pro}\label{Prop-Main k}\textit{
Let $M^n$ be a complete Riemannian manifold with $\Ric \geq 0$. For
$k\in \N$, there exists a constant $\delta_k(n) > 0$ such that if
$\alpha_M \geq 1 - \delta_k(n)$, then $\pi_k(M^n) = 0$.}
\end{Pro}

\begin{proof}
Choose some $c>1$ and set $\delta_k(n) = 1 - \beta(k,c,n)$. The
conclusion then follows from Lemma \ref{Lemma-Main}.
\end{proof}

Thus, we recover Perelman's result \cite{P2}:

\begin{Lem}\label{Prop-Perelman Thm 2}{\em [\cite{P2}, Theorem 2]}.
Let $M^n$ be a complete Riemannian manifold with $\Ric \geq 0$.
There exists a constant $\delta_n >0$ such that if $\alpha_M \geq 1
- \delta_n$, then $M^n$ is contractible.
\end{Lem}

\begin{proof}
Choose some $c>1$ and set $\delta_n = 1 - \max_{k=1,..,n}
\beta(k,c,n)$. Then Lemma \ref{Lemma-Main} implies $\pi_k(M^n) =0$
for all positive values $k$. Hence, $M^n$ is contractible by the
Whitehead Theorem \cite {W}.
\end{proof}

{\bf Remark.} In the appendix we use the expression for
$\beta(k,c,n)$ from Definition \ref{Def-beta volume growth} to find
the `best' value (depending only on $k$ and $n$) of $\alpha_M$ which
guarantees that $\pi_k(M^n) =0$. This is the lower bound for
$\alpha_M$ as stated in Theorem \ref{Theorem-Main}.

We now prove Theorem \ref{Theorem-Main}.
\begin{proof}
Let
$$\alpha(k,n) = \inf_{c \in (1, \infty)} \beta(k,c,n).$$
By assumption, $\alpha_M > \alpha(k,n)$ and thus there exists
$c_0>1$ such that $\alpha_M \geq \beta(k, c_0, n)$. The result
follows by applying Main Lemma($k$) with $c=c_0$. In the appendix we
compute values of $\alpha(k,n)$.
\end{proof}

\section{Appendix}\label{Section-Appendix}
The constants $\emph{C}_{k,n}$ explicitly determine the function $h_{k,n}(x)$ defined in section \ref{Section-General Argument}. In this appendix, we show that the constants $\emph{C}_{k,n}$ as defined are optimal and use the definition of $h_{k,n}$ to compute explicit values for $\alpha(k,n)$ as stated in Theorem \ref{Theorem-Main}.

\subsection{Optimal Constants}
Recall Definition \ref{Def-C(k,n) constants} of $\emph{C}_{k,n}(i)$:
\begin{equation}
\emph{C}_{k,n}(i) = (16k)^{n-1}(1+10\emph{C}_{k,n}(i-1))^n + 3 +
10\emph{C}_{k,n}(i-1), \qquad i \geq 1
\end{equation}
and $\emph{C}_{k,n}(0)=1$. Denote $\emph{C}_{k,n} =
\emph{C}_{k,n}(k)$

The constants $\emph{C}_{k,n}$ grow large very quickly. Preliminary values for $\emph{C}_{k,n}$ where $1\leq k \leq 3$ and $1 \leq n \leq 8$ are listed in Table \ref{Table of constants}.

\begin{table}[h]
\caption{\small {Table of $\emph{C}_{k,n}$ values for $1 \leq k \leq 3, 1 \leq n \leq 10$}}
\label{Table of constants}
\begin{center}
\begin{tabular}{|c||c|c|c|}
  \hline
   & $k=1$ & $k=2$ & $k=3$ \\  \hline \cline{1-4}
  $n=1$ & $24$ & - & - \\ \hline
  $n=2$ & $384$ & $1.89 \times 10^{8}$ & - \\ \hline
  $n=3$ & $6144$ & $1.52 \times 10^{17}$ & $1.36 \times 10^{60}$ \\ \hline
  $n=4$ & $98304$ & $1.25 \times 10^{29}$ & $1.00 \times 10^{133}$ \\ \hline
  $n=5$ & $1.57 \times 10^{6}$ & $1.06 \times 10^{44}$ & $9.53 \times 10^{248}$ \\ \hline
  $n=6$ & $2.51 \times 10^{7}$ & $9.15 \times 10^{61}$ & $1.43 \times 10^{418}$ \\ \hline
  $n=7$ & $4.03 \times 10^{8}$ & $8.10 \times 10^{82}$ & $4.12 \times 10^{650}$ \\ \hline
  $n=8$ & $6.44 \times 10^{9}$ & $7.35 \times 10^{106}$ & $2.80 \times 10^{956}$ \\ \hline
  $n=9$ & $1.03 \times 10^{11}$ & $6.82 \times 10^{133}$ & $5.50 \times 10^{1345}$ \\ \hline
  $n=10$ & $1.65 \times 10^{12}$ & $6.49 \times 10^{163}$ & $3.81 \times 10^{1828}$ \\ \hline
\end{tabular}
\end{center}
\end{table}

\begin{Lem} \label{Lemma-d_i and b_i}\textit{
If $d_i = \emph{C}_{k,n}(i) d_0$ and $b_i =
[16k(1+10\emph{C}_{k,n}(i))]^{-(n-1)},$ then
\begin{equation}\label{Eq1}
d_0 + 10d_i = b_i(d_{i+1} -3d_0 -10d_i),
\end{equation}
and
\begin{equation}\label{Eq3}
8b_i^{\frac{1}{n-1}}(d_0 + 10d_i) = \frac{d_0}{2k},
\end{equation}
for $i=0,1,..,k$. Furthermore, (\ref{Ineq2}) and (\ref{Ineq4}) hold
as well.}
\end{Lem}

\begin{proof} The proofs of (\ref{Eq1}) and (\ref{Eq3}) are by induction in $i$. When $i=0$ the conclusion holds.
Assume the conclusion holds for $i < k$. It remains to verify the
conclusion for $i+1$. Note that
\begin{eqnarray*}
&& \hspace{-1.6cm} b_{i+1}(d_{i+2}-3d_0-10d_{i+1}) \\
&~~=& [16k(1+10\emph{C}_{k,n}(i+1))]^{n-1}(\emph{C}_{k,n}(i+2)d_0 - 3d_0 -
10\emph{C}_{k,n}(i+1)d_0)\\
&~~=&(1+10\emph{C}_{k,n}(i+1))d_0 = d_0 + 10d_{i+1}.
\end{eqnarray*}

Similarly, for the second equation we get
$$8b_{i+1}^{\frac{1}{n-1}}(d_0+10d_{i+1})=8[16k(1+10\emph{C}_{k,n}(i+1))]^{-1}(1+10\emph{C}_{k,n}(i+1))d_0=\frac{d_0}{2k}.$$

To verify that (\ref{Ineq4}) holds, note that $d_k = C_{k,n}(k)d_0 =
C_{k,n}d_0$ and, by the definition of $h_{k,n}(d_0)$ [Definition
\ref{Def-six* function}], we have exactly (\ref{Ineq4}).

Lastly, both (\ref{Ineq4}) and (\ref{Eq1}) imply (\ref{Ineq2}). Note
that, setting $h_{k,n}(d_0)=c$, (\ref{Ineq4}) implies
\begin{eqnarray}
d_k &\leq& 10^{-k-2} (1+d_0/2k)^{-k}(1-c^{-1})\\
&=& 10^{-k-2} (1+d_0/2k)^{-k} 1/c (c-1)\\
&\leq& c-1,
\end{eqnarray}
where the last inequality follows because $10^{-k-2} < 1$,
$(1+d_0/2k)^{-k} < 1$, and $1/c < 1$.

Therefore, since $1 \leq i < k$,
\begin{eqnarray}
b_i(c-1 + d_0(2-i/k)) &\geq& b_i(d_k + d_0(2-i/k))\\
&\geq& b_i(d_k + d_0)\\
&\geq& b_i \cdot d_k\\
&\geq& b_i \cdot d_{i+1}\\
&\geq& b_i (d_{i+1} - 3d_0 - 10d_i)\\
&=& d_0 + 10d_i,
\end{eqnarray}
where the last equality follows from (\ref{Eq1}). Thus,
(\ref{Ineq2}) holds and this completes the proof.
\end{proof}

{\bf Remark.} So we see that the constants $\emph{C}_{k,n}(i)$ suffice for the proof of Theorem \ref{Theorem-Main}.
Next we show that these constants provide the optimal choice.

\begin{Lem}\label{Lemma-Appendix1}\textit{
If (\ref{Ineq1}) and (\ref{Ineq3}) hold for all $i\geq0$, then
$$d_i \geq \emph{C}_{k,n}(i)$$ and $$b_i \leq
\frac{1}{[16k(1+10\emph{C}_{k,n}(i))]^{n-1}}.$$}
\end{Lem}

\begin{proof}
The proof is by induction on $i$. When $i=0$ the conclusion holds.
From (\ref{Ineq1}) and assuming the conclusion holds for $i$, we
have
\begin{eqnarray}
d_{i+1} &\geq& \frac{1}{b_i}(d_0 + 10d_i) + 3d_0 + 10d_i\\
&\geq& [(16k)^{n-1}(1+10\emph{C}_{k,n}(i))^n + 3 + 10
\emph{C}_{k,n}(i)]d_0\\
&=& \emph{C}_{k,n}(i+1)d_0.
\end{eqnarray}
Using this lower bound for $d_{i+1}$ and (\ref{Ineq3}, we get
\begin{eqnarray}
b_{i+1} &\leq& \left(\frac{d_0}{2k} \frac{1}{d_0 + 10d_{i+1}}
\right)^{n-1}\\
&=& \left(\frac{1}{16k(1+10\emph{C}_{k,n}(i+1))} \right)^{n-1}.
\end{eqnarray}
This completes the proof.
\end{proof}

\subsection{Computing $\alpha(k,n)$ values}\label{section-computing alpha(k,n) values}

The term $\beta(k,c,n)$ denotes the minimal volume growth necessary
to guarantee that any continuous map $f: S^k \to B_p(R)$ has a
continuous extension $g: D^{k+1} \to B_p(cR)$ (see Definition
\ref{Def-beta volume growth}). Recall that the expression for
$\beta(k,c,n)$ is iteratively defined.

By definition,
\begin{eqnarray}
\beta(k,c,n) &= \max \Bigg\{&1-\gamma\left(c,h_{k,n}^{-1}\left(c\right),n\right);\\
&&\beta\left(j, 1+\frac{h_{k,n}^{-1}\left(c\right)}{2k},n\right),
j=1,..,k-1 \Bigg\}.
\end{eqnarray}

Ultimately we are not concerned with the location of the homotopy
map. Thus we have a certain amount of freedom when choosing which
$c$ value to take. To determine the optimal bound on volume growth
guaranteeing $\pi_k(M^n)=0$, it is necessary to choose the `best'
value of $c$ for $\beta(k,c,n)$; that is, the $c$ which makes
$\beta(k,c,n)$ the smallest. Set $\alpha(k,n) = \inf_{c >1}
\beta(k,c,n)$.

In order to compute explicit values for $\alpha(k,n)$, we must
successively simplify the components of $\beta(k,c,n)$. Ultimately,
because of its iterative definition, it is possible to express
$\beta(k,c,n)$ as the maximum of a collection of $\gamma$ terms.
Using the definition of $\gamma(c,\epsilon, n)$, we can then compute
specific values for $\alpha(k,n)$. Here we describe in detail the
method to compute $\alpha(k,n)$ and compile a table of these values
for $k=1,2,3$ and $n=1,...,10$.

To begin, we have
\begin{equation}\label{beta1}
\beta(1,c,n) = 1-\gamma\left(c,h_{1,n}^{-1}\left(c\right),n\right).
\end{equation}
By definition, when $k=2$
\begin{eqnarray}
\beta(2,c,n) &= \max\Bigg\{&1-\gamma\left(c,h_{2,n}^{-1}\left(c\right),n\right),\\
&&\beta\left(1, 1+\frac{h_{2,n}^{-1}\left(c\right)}{4},n\right)\Bigg\}.
\end{eqnarray}
To evaluate this expression for $\beta(2,c,n)$, simplify the
$\beta\left(1,1+\frac{h_{2,n}^{-1}\left(c\right)}{4},n\right)$ term
by setting $c=1+\frac{h_{2,n}^{-1}\left(c\right)}{4}$ and applying
(\ref{beta1}). Therefore,

\begin{eqnarray}\label{beta2}
\hspace{-.4cm}
\beta(2,c,n) &\hspace{-.2cm} = \hspace{-.1cm} \max \Bigg\{& \hspace{-.3cm} 1-\gamma\left(c,h_{2,n}^{-1}\left(c\right),n\right),\\
&&
\hspace{-.3cm} 1-\gamma\left(1+\frac{h_{2,n}^{-1}\left(c\right)}{4},h_{1,n}^{-1}\left(1+\frac{h_{2,n}^{-1}\left(c\right)}{4}\right),n\right)\Bigg\}.
\end{eqnarray}

Similarly, to evaluate $\beta(3,c,n)$ we have, by definition,
\begin{eqnarray}
\beta(3,c,n) &=\max\Bigg\{& \hspace{-.2cm} 1-\gamma\left(c,h_{3,n}^{-1}\left(c\right),n\right);\\
&& \hspace{-.2cm} \beta\left(j, 1+\frac{h_{3,n}^{-1}\left(c\right)}{6},n\right),
j=1,2\Bigg\}\\
&= \max\Bigg\{& \hspace{-.2cm} 1-\gamma\left(c,h_{3,n}^{-1}\left(c\right),n\right),\\
&& \hspace{-.2cm} \beta\left(1, 1+\frac{h_{3,n}^{-1}\left(c\right)}{6},n\right),\\
&& \hspace{-.2cm} \beta\left(2, 1+\frac{h_{3,n}^{-1}\left(c\right)}{6},n\right)\Bigg\}.
\end{eqnarray}
Substituting $\beta\left(1,1+\frac{h_{3,n}^{-1}\left(c\right)}{6},n\right)$ with the expression obtained by setting\\
$c=1+\frac{h_{3,n}^{-1}\left(c\right)}{6}$ and evaluating (\ref{beta1}) yields
\begin{eqnarray}
\beta(3,c,n) &\hspace{-.2cm} = \max\Bigg\{& \hspace{-.3cm} 1-\gamma\left(c,h_{3,n}^{-1}\left(c\right),n\right),\\
&& \hspace{-.3cm} 1-\gamma\left(1+\frac{h_{3,n}^{-1}\left(c\right)}{6},h_{1,n}^{-1}\left(1+\frac{h_{3,n}^{-1}\left(c\right)}{6}\right),n\right),\\
&& \hspace{-.3cm} \beta\left(2, 1+\frac{h_{3,n}^{-1}\left(c\right)}{6},n\right)\Bigg\}.
\end{eqnarray}

Finally, apply (\ref{beta2}) with
$c=1+\frac{h_{3,n}^{-1}\left(c\right)}{6}$ to simplify the
remaining\\
$\beta\left(2,1+\frac{h_{3,n}^{-1}\left(c\right)}{6},n\right)$ term.
We get

\begin{eqnarray*}
\noindent
&&\hspace{-2.4cm} \beta(3,c,n) \\
&= \max \hspace{-.1cm}\Bigg\{& \hspace{-.3cm} 1-\gamma\left(c,h_{3,n}^{-1}\left(c\right),n\right),\\
&& \hspace{-.4cm} 1-\gamma\left(1+\frac{h_{3,n}^{-1}\left(c\right)}{6},h_{1,n}^{-1}\left(1+\frac{h_{3,n}^{-1}\left(c\right)}{6}\right),n\right),\\
&& \hspace{-.4cm} 1-\gamma\left(1+\frac{h_{3,n}^{-1}\left(c\right)}{6},h_{2,n}^{-1}\left(1+\frac{h_{3,n}^{-1}\left(c\right)}{6}\right),n\right),\\
&&
\hspace{-.4cm} 1-\gamma\left(1+\frac{h_{2,n}^{-1}\left(1+\frac{h_{3,n}^{-1}\left(c\right)}{6}\right)}{4},h_{1,n}^{-1}\left(1+\frac{h_{2,n}^{-1}\left(1+\frac{h_{3,n}^{-1}\left(c\right)}{6}\right)}{4}\right),n\right)\Bigg\}.
\end{eqnarray*}

Because of the successive nesting, when completely expanded, the
expression $\beta(k,c,n)$ can be written as the maximum of $2^{k-1}$
terms of the form $1 - \gamma(.,.,n)$. However, given the nature of
the functions $h_{k,n}(x)$ and the behavior of $\gamma(c,
h_{k,n}^{-1}(c), n)$ when $1 < c < 2$, the maximum of this
collection of $1 - \gamma$ terms is determined by the the maximum of
the leading $1 - \gamma(c, h_{k,n}^{-1}(c), n)$ term and the final
$1 - \gamma$ term containing the most iterations. That is to say,
for all $k$ and $n$,

\begin{eqnarray*}
&\hspace{-2cm} \beta(k,c,n)& \\
&& \hspace{-1.8cm}
= \max  \Bigg\{  1-\gamma\left(c,h_{k,n}^{-1}\left(c\right),n\right); \beta\left(j, 1+\frac{h_{k,n}^{-1}\left(c\right)}{2k},n\right),
j=1,..,k-1 \Bigg\} \\
&\hspace{-.4cm}= \max \Bigg\{ & \hspace{-.4cm} 1-\gamma\left(c,h_{k,n}^{-1}\left(c\right),n\right),\\
&& \hspace{-1cm} 1-\gamma\left( 1 + \dots \frac{h_{k-1,n}^{-1}\left(1 + \frac{h_{k,n}^{-1}(c)}{2k}  \right)}{2(k-1)}, h_{1,n}^{-1}\left( 1 + \dots \frac{h_{k-1,n}^{-1}\left(1 + \frac{h_{k,n}^{-1}(c)}{2k}  \right)}{2(k-1)} \right), n \right) \Bigg\};
\end{eqnarray*}
\noindent
which in turn can be written as
\vspace{.1cm}
\begin{eqnarray*}
&& \hspace{-1cm} \beta(k,c,n) = \max \Bigg\{ 1-\gamma\left(c,h_{k,n}^{-1}\left(c\right),n\right),\\
&& \hspace{.6cm} 1-\gamma\left( 1 + \dots \frac{h_{k-1,n}^{-1}\left(1 + \frac{h_{k,n}^{-1}(c)}{2k}  \right)}{2(k-1)}, h_{1,n}^{-1}\left( 1 + \dots \frac{h_{k-1,n}^{-1}\left(1 + \frac{h_{k,n}^{-1}(c)}{2k}  \right)}{2(k-1)} \right), n \right) \Bigg\}\\
&& \hspace{.6cm} = 1 - \min \Bigg\{ \gamma\left(c,h_{k,n}^{-1}\left(c\right),n\right),\\
&& \hspace{.6cm} \gamma\left( 1 + \dots \frac{h_{k-1,n}^{-1}\left(1 + \frac{h_{k,n}^{-1}(c)}{2k}  \right)}{2(k-1)}, h_{1,n}^{-1}\left( 1 + \dots \frac{h_{k-1,n}^{-1}\left(1 + \frac{h_{k,n}^{-1}(c)}{2k}  \right)}{2(k-1)} \right), n \right) \Bigg\}.
\end{eqnarray*}
\vspace{.1cm} \indent Recall that the constants $\delta_{k,n}$
[Definition \ref{Def-six* function}] represent the location of the
vertical asymptote $x = \delta_{k,n}$ of the function $h_{k,n}(x)$.
Therefore, the function $h_{k,n}^{-1}$ is bounded above by the
constant $\delta_{k,n}$; that is, $h_{k,n}^{-1}(c) < \delta_{k,n}$
for all $c>1$. In Table \ref{Table of deltas}, we list values of
$\delta_{k,n}$ for $1 \leq k \leq 3$ and $1 \leq n \leq 10$.

\begin{table}[h]
\caption{{Table of $\delta_{k,n}$ values for $1 \leq k \leq 3, 1 \leq n \leq 10$}}
\label{Table of deltas}
\begin{center}
\begin{tabular}{|c||c|c|c|}
  \hline
   & $k=1$ & $k=2$ & $k=3$ \\  \hline \cline{1-4}
  $n=1$ & $4.17 \times 10^{-5}$ & - & - \\ \hline
  $n=2$ & $2.60 \times 10^{-6}$ & $5.29 \times 10^{-13}$ & - \\ \hline
  $n=3$ & $1.63 \times 10^{-7}$ & $6.58 \times 10^{-22}$ & $7.34 \times 10^{-66}$ \\ \hline
  $n=4$ & $1.02 \times 10^{-8}$ & $7.98 \times 10^{-34}$ & $9.96 \times 10^{-139}$ \\ \hline
  $n=5$ & $6.36 \times 10^{-10}$ & $9.45 \times 10^{-49}$ & $1.05 \times 10^{-254}$ \\ \hline
  $n=6$ & $3.97 \times 10^{-11}$ & $1.09 \times 10^{-66}$ & $7.01 \times 10^{-424}$ \\ \hline
  $n=7$ & $2.48 \times 10^{-12}$ & $1.23 \times 10^{-87}$ & $2.43 \times 10^{-656}$ \\ \hline
  $n=8$ & $1.55 \times 10^{-13}$ & $1.36 \times 10^{-111}$ & $3.57 \times 10^{-962}$ \\ \hline
  $n=9$ & $9.70 \times 10^{-15}$ & $1.47 \times 10^{-138}$ & $1.81 \times 10^{-1351}$ \\ \hline
  $n=10$ & $6.06 \times 10^{-16}$ & $1.54 \times 10^{-168}$ & $2.62 \times 10^{-1834}$ \\
  \hline
\end{tabular}
\end{center}
\end{table}

Fixing $k$ and $n$, the function $\gamma(c, h_{k,n}^{-1}(c), n)$ is
increasing as a function of $c$ when $1<c<2$. Further, we have that
for all $c > 1$

\begin{eqnarray}
h_{k,n}^{-1}(c) &<& \delta_{k,n}\\
h_{k,n}^{-1}(c)/2k &<& \delta_{k,n}/2k\\
1 + h_{k,n}^{-1}(c)/2k &<& 1 + \delta_{k,n}/2k << 2.
\end{eqnarray}

Define $\epsilon_{k,n}$ as

\begin{eqnarray*}
&& \hspace{-.5cm} \epsilon_{k,n} = \lim_{c \to \infty} \gamma\left( 1 + \dots \frac{h_{k-1,n}^{-1}\left(1 + \frac{h_{k,n}^{-1}(c)}{2k}  \right)}{2(k-1)}, h_{1,n}^{-1}\left( 1 + \dots \frac{h_{k-1,n}^{-1}\left(1 + \frac{h_{k,n}^{-1}(c)}{2k}  \right)}{2(k-1)} \right), n \right)
\end{eqnarray*}
\begin{eqnarray*}
&=& \gamma\left( 1 + \dots \frac{h_{k-1,n}^{-1}\left(1 + \frac{\delta_{k,n}}{2k}  \right)}{2(k-1)}, h_{1,n}^{-1}\left( 1 + \dots \frac{h_{k-1,n}^{-1}\left(1 + \frac{\delta_{k,n}}{2k}  \right)}{2(k-1)} \right), n \right)\\
&=& \left[1 + \left(\frac{1 + \dots \frac{h_{k-1,n}^{-1}\left(1 + \frac{\delta_{k,n}}{2k}  \right)}{2(k-1)}}{h_{1,n}^{-1}\left( 1 + \dots \frac{h_{k-1,n}^{-1}\left(1 + \frac{\delta_{k,n}}{2k}  \right)}{2(k-1)} \right)}\right)^n \right]^{-1}.
\end{eqnarray*}

With this simplification, it is possible to explicitly compute the
values of $\epsilon_{k,n}$. Table \ref{Table of epsilons} below
lists values of $\epsilon_{k,n}$ for $k=1,2,3$ and $n=1,...,10$.
These values were computing using Mathematica 6.0 and the source
code for these computations as well as additional exposition is available in \cite{Mu}.

\begin{table}[h]
\caption{\small{Table of $\epsilon_{k,n}$ values for $1 \leq k \leq 3, 1 \leq n \leq 10$}}
\label{Table of epsilons}
\begin{center}
\begin{tabular}{|c||c|c|c|}
  \hline
   & $k=1$ & $k=2$ & $k=3$ \\  \hline \cline{1-4}
  $n=1$ & $1.04 \times 10^{-5}$ & - & - \\ \hline
  $n=2$ & $4.24 \times 10^{-13}$ & $1.89 \times 10^{-37}$ & - \\ \hline
  $n=3$ & $6.74 \times 10^{-23}$ & $1.92 \times 10^{-86}$ & $3.52 \times 10^{-284}$ \\ \hline
  $n=4$ & $4.18 \times 10^{-35}$ & $1.70 \times 10^{-167}$ & $1.29 \times 10^{-722}$ \\ \hline
  $n=5$ & $1.01 \times 10^{-49}$ & $7.64 \times 10^{-290}$ & $1.25 \times 10^{-1563}$ \\ \hline
  $n=6$ & $9.61 \times 10^{-67}$ & $1.64 \times 10^{-462}$ & $4.16 \times 10^{-3006}$ \\ \hline
  $n=7$ & $3.56 \times 10^{-86}$ & $1.55 \times 10^{-694}$ & $2.75 \times 10^{-5289}$ \\ \hline
  $n=8$ & $5.14 \times 10^{-108}$ & $6.06 \times 10^{-995}$ & $9.42 \times 10^{-8693}$ \\ \hline
  $n=9$ & $2.90 \times 10^{-132}$ & $9.08 \times 10^{-1373}$ & $1.94 \times 10^{-13536}$ \\ \hline
  $n=10$ & $6.41 \times 10^{-159}$ & $4.87 \times 10^{-1837}$ & $1.24 \times 10^{-20180}$ \\ \hline
\end{tabular}
\end{center}
\end{table}

\newpage

The value $\alpha(k,n)$, as described in Theorem \ref{Theorem-Main}, represents the optimal lower bound for the volume growth guaranteeing  $\pi_k(M^n) =0$. We can then set $\alpha(k,n) = 1 - \epsilon_{k,n}$. Table \ref{Table of alphas} contains the values of $\alpha(k,n)$ for $k=1,2,3$ and $n=1,...,10$.

\begin{table}[!h]
\caption{\small{Table of $\alpha(k,n)$ values for $1 \leq k \leq 3, 1 \leq n \leq 10$}}
\label{Table of alphas}
\begin{center}
\begin{tabular}{|c||c|c|c|}
  \hline
   & $k=1$ & $k=2$ & $k=3$ \\ \hline \cline{1-4}
  $n=1$ & $1 - 1.04 \times 10^{-5}$ & - & - \\ \hline
  $n=2$ & $1 - 4.24 \times 10^{-13}$ & $1 - 1.89 \times 10^{-37}$ & - \\ \hline
  $n=3$ & $1 - 6.74 \times 10^{-23}$ & $1 - 1.92 \times 10^{-86}$ & $1 - 3.52 \times 10^{-284}$ \\ \hline
  $n=4$ & $1 - 4.18 \times 10^{-35}$ & $1 - 1.70 \times 10^{-167}$ & $1 - 1.29 \times 10^{-722}$ \\ \hline
  $n=5$ & $1 - 1.01 \times 10^{-49}$ & $1 - 7.64 \times 10^{-290}$ & $1 - 1.25 \times 10^{-1563}$ \\ \hline
  $n=6$ & $1 - 9.61 \times 10^{-67}$ & $1 - 1.64 \times 10^{-462}$ & $1 - 4.16 \times 10^{-3006}$ \\ \hline
  $n=7$ & $1 - 3.56 \times 10^{-86}$ & $1 - 1.55 \times 10^{-694}$ & $1 - 2.75 \times 10^{-5289}$ \\ \hline
  $n=8$ & $1 - 5.14 \times 10^{-108}$ & $1 - 6.06 \times 10^{-995}$ & $1 - 9.42 \times 10^{-8693}$ \\ \hline
  $n=9$ & $1 - 2.90 \times 10^{-132}$ & $1 - 9.08 \times 10^{-1373}$ & $1 - 1.94 \times 10^{-13536}$ \\ \hline
  $n=10$ & $1 - 6.41 \times 10^{-159}$ & $1 - 4.87 \times 10^{-1837}$ & $1 - 1.24 \times 10^{-20180}$ \\ \hline
\end{tabular}
\end{center}
\end{table}

In general, $\alpha(1, n) = 1 - \left[1 +
\frac{2}{h^{-1}_{1, n}(2)} \right]^{-1}$; and for $k > 1$, we have
\begin{eqnarray}\label{alphas}
\alpha(k,n) &=& 1 - \epsilon_{k,n}\\
& = & 1 - \left[1 + \left(\frac{1 + \dots \frac{h_{k-1,n}^{-1}\left(1 + \frac{\delta_{k,n}}{2k}  \right)}{2(k-1)}}{h_{1,n}^{-1}\left( 1 + \dots \frac{h_{k-1,n}^{-1}\left(1 + \frac{\delta_{k,n}}{2k}  \right)}{2(k-1)} \right)}\right)^n \right]^{-1}.
\end{eqnarray}
These are the bounds are the best that can be achieved via Perelman's method \cite{P2}.

Combining this information with previous results of Anderson \cite{A}, Li \cite{Li}, Cohn-Vossen \cite{Cohn} and Zhu \cite{Z1} we can refine the table above.

\begin{table}[h]
\caption{\small{Table of revised $\alpha(k,n)$ values for $1 \leq k \leq 3, 1 \leq n \leq 10$}}
\begin{center}
\begin{tabular}{|c||c|c|c|}
  \hline
   & $k=1$ & $k=2$ & $k=3$ \\ \hline \cline{1-4}
  $n=1$ & $-$ & - & - \\ \hline
  $n=2$ & $0$ & $0$ & - \\ \hline
  $n=3$ & $0$ & $0$ & $0$ \\ \hline
  $n=4$ & $1/2$ & $1 - 1.70 \times 10^{-167}$ & $1 - 1.29 \times 10^{-722}$ \\ \hline
  $n=5$ & $1/2$ & $1 - 7.64 \times 10^{-290}$ & $1 - 1.25 \times 10^{-1563}$ \\ \hline
  $n=6$ & $1/2$ & $1 - 1.64 \times 10^{-462}$ & $1 - 4.16 \times 10^{-3006}$ \\ \hline
  $n=7$ & $1/2$ & $1 - 1.55 \times 10^{-694}$ & $1 - 2.75 \times 10^{-5289}$ \\ \hline
  $n=8$ & $1/2$ & $1 - 6.06 \times 10^{-995}$ & $1 - 9.42 \times 10^{-8693}$ \\ \hline
  $n=9$ & $1/2$ & $1 - 9.08 \times 10^{-1373}$ & $1 - 1.94 \times 10^{-13536}$ \\ \hline
  $n=10$ & $1/2$ & $1 - 4.87 \times 10^{-1837}$ & $1 - 1.24 \times 10^{-20180}$ \\ \hline
\end{tabular}
\end{center}
\end{table}

\singlespacing
\vspace{5mm}\noindent
Michael Munn \\
New York City College of Technology, CUNY \\
Brooklyn, NY\\
e-mail: mikemunn@gmail.com

\end{document}